\documentclass[11pt]{article}
\usepackage{amsmath, amsthm, amssymb, bm, booktabs, mathtools, multirow, graphicx, subfigure, xcolor}
\usepackage[margin=1in]{geometry}
\usepackage{listings}
\usepackage{natbib}
\usepackage{natbib}
\setcitestyle{authoryear,open={(},close={)}} 
\usepackage{hyperref}

\usepackage{ifthen}
\usepackage{soul}

\newtheorem{definition}{Definition}[section]
\newtheorem{theorem}{Theorem}[section]

\newtheorem{corollary}{Corollary}[section]
\newtheorem{proposition}{Proposition}[section]
\newtheorem{example}{Example}[section]
\newtheorem{remark}{Remark}[section]

\numberwithin{equation}{section}

\newcommand*\D{\,\mathrm{d}}

\newcommand{\replace}[2]{{#2}}

\makeatletter
\newcommand*\bigcdot{\mathpalette\bigcdot@{.5}}
\newcommand*\bigcdot@[2]{\mathbin{\vcenter{\hbox{\scalebox{#2}{$\m@th#1\bullet$}}}}}
\makeatother

\providecommand{\keywords}[1]
{
  \small	
  \textbf{\textit{Keywords---}} #1
}

\title{Partial Conditioning for Inference of Many-Normal-Means with H\"older Constraints}
\author{Jiasen Yang, Xiao Wang, Chuanhai Liu\\
Department of Statistics\\
Purdue University}
\date{}

\begin{document}

\maketitle
\vspace{-2.5in}
\noindent
\textit{International Journal of Approximate Reasoning} \textbf{159}\\
\url{https://doi.org/10.1016/j.ijar.2023.108946}
\vspace{2in}
\begin{abstract}


Inferential models have been proposed for valid and efficient prior-free
probabilistic inference.  As it gradually gained popularity, this theory is subject
to further developments for practically challenging problems.
This paper considers the many-normal-means problem with the means
constrained to be in the neighborhood of each other\replace{}{, formally represented by a H\"older space.}
A new method, called partial conditioning, is proposed
to generate valid and efficient marginal inference about the individual means.
It is shown that the method outperforms both a fiducial-counterpart in terms of validity
and a conservative-counterpart in terms of efficiency.
\replace{}{We conclude the paper by remarking that a general theory of partial conditioning for inferential models deserves future development.}

\end{abstract}




\section{Introduction}
\label{s:Introduction}

The framework of inferential models (IMs) has been developed to provide
valid and efficient prior-free probabilistic inference.
In its simplest form \citep{martin2013inferential}, it is obtained
by modifying R.~A.~Fisher's fiducial
argument, which suffers from the mathematical difficulties discussed, for example, in \cite{liu2015frameworks}.
It has also been extended for efficient inference by combining information
\citep{martin2015conditional}
and marginalization \citep{martin2015marginal}.
For certain practically challenging problems,
however,
this theory is subject to further developments,
similar to other existing schools of thought.
One noticeable class of such challenging problems 
is inference on constrained parameters, which is challenging for all existing inferential methods (see, {\it e.g.},
\cite{leaf2012inference} and references therein).

In this paper, we consider the many-normal-means problem with the means
constrained to be in the neighborhood of each other.
For both conceptual and representational simplicity,
we describe the problem as a special case of the familiar classic
formulation of a nonparametric regression model \citep{lepskii1991problem}.
More precisely, to motivate the problem, we consider the following model
\begin{equation}\label{eq:model}
Y_i  = \vartheta_0(t_i) + \sigma U_i, \quad i=1, \cdots, n,
\end{equation}
where the $t_i$ denotes the design points, the $Y_i$ denotes the responses, and $U_i\stackrel{i.i.d.}{\sim} \mathcal{N}(0,1)$ denote the random error terms. Here, $\vartheta_0$ is the unknown regression function that is assumed to reside in the H\"older space
\begin{equation}\label{eq:Holder}
\Theta_{M,\gamma} = \left\{\vartheta: [0, 1]\rightarrow \mathbb R\ \Big|\ |\vartheta(t) - \vartheta(s)|\le M |t-s|^{\gamma},\ \forall t,s\in[0,1] \right\},
\end{equation}
with the H\"older exponent $0<\gamma\le 1$.
We assume $M$, $\gamma$, and $\sigma$ to be known constants and, without loss of generality, we set $\sigma=1$.
In subsequent discussions, we suppress the subscripts and denote the H\"older space $\Theta_{M,\gamma}$ simply by $\Theta$, unless otherwise noted. 
\replace{}{The motivation for this problem stems from the challenging setting in nonparametric regression, which we briefly discuss in the following remark.}

\begin{remark}
\replace{}{
In nonparametric regression, the objective is to estimate $f$ from a nonparametric class of functions ${\cal F}$, under the assumption that $f$ belongs to this class. For instance, ${\cal F}$ can be the set of all continuous functions on $[0, 1]$. In this paper, we focus on the H\"older class $\Theta_{M,\gamma}$. 
 A function in $\Theta_{M,\gamma}$ with $\gamma>1$ is constant, while a function satisfying $\gamma=1$ meets a Lipschitz condition with a Lipschitz constant of~$M$. Thus, the space of Lipschitz continuous functions is a subset of $\Theta_{M,\gamma}$. It is very challenging to construct confidence intervals for $f$ either locally or globally. For example, consider the Nadaraya-Watson estimator $\hat \vartheta_n$ with kernel-bandwidth $h$. A Taylor expansion shows that it has bounded bias and variance
\[
|\mathbb E(\hat \vartheta_n(x)) - \vartheta_0(x)| \le C_1 h^\gamma, \quad
\mathrm{Var}(\hat \vartheta_n(x)) = {C_2\over nh}\,,
\]
for all $x\in [0, 1]$ with some positive constants $C_1$ and $C_2$. By taking $h = O(n^{-1/(2\gamma+1)})$, the $L_2$ error of $\hat \vartheta_n$ achieves the optimal convergence rate (cf.~Section 1.6.1 of \cite{Tsybakov})
\[
\mathbb E\|\hat \vartheta_n - \vartheta_0\|_{L_2}^2 = O(n^{-2\gamma/(2\gamma+1)})\,
\]
or
\begin{equation}\label{eq:asym-rate}
\sqrt{\mathbb E\|\hat \vartheta_n - \vartheta_0\|_{L_2}^2} = O(n^{-\gamma/(2\gamma+1)})\,.
\end{equation}
However, under the optimal $h$, the asymptotic bias has the same order as the asymptotic standard deviation, which is of order $O(n^{-\gamma/(2\gamma+1)})$. To eliminate the bias, a typical strategy is to opt for a smaller order $h$, which leads to a suboptimal asymptotic confidence interval in terms of convergence rate.
}
\end{remark}

Our objective here in this paper is to perform statistical inference on the unknown function $\vartheta_0$ within the IM framework.
We propose a new method, called \emph{partial conditioning}, to generate valid and efficient marginal inference about the individual means.
We show that the method outperforms both a fiducial-counterpart in term of validity
and a conservative-counterpart in terms of efficiency. \replace{}{A simple simulation-based study shows that our empirical asymptotic convergence rate of plausibility intervals for the case of $\gamma=\frac{1}{2}$ is about $O(n^{-0.267})$, slightly better than or at least close to that in \eqref{eq:asym-rate}, which in this case is $O(n^{-0.250})$.  While one of these two convergence rates is on point estimation and the other on interval length, such a comparison is arguably meaningful because IM plausibility intervals are valid in terms of frequency calibration.}

It should be noted that IM approaches to the general nonparametric problem itself deserve in-depth investigations. This explains why our focus here is on a simple case. Our goal is to introduce an innovative idea of conditioning in the IM framework to combine information when multiple parameters entangle with each other via known constraints. \replace{}{As noted in the above remark on the difficulty of nonparametric regression, no existing methods can produce valid inference, even asymptotically. We believe our proposed method of partial conditioning makes a promising step in developing satisfactory solutions to nonparametric regression.}


The rest of the paper is arranged as follows.
Section \ref{s:BasicAndConditionalIMs} reviews basic IMs \citep{martin2013inferential}
and conditional IMs \citep{martin2015conditional}
with the simple cases of $n=1$ and $n=2$ with $M=0$ in \eqref{eq:model} and 
\eqref{eq:Holder}.
Section \ref{s:PartialConditioning}
addresses the difficulties in the case of $n=2$ and $M\neq 0$
and introduces the proposed method of partial conditioning,
with the general $n$ case considered in
Section \ref{s:GeneralCase}, where illustrative numerical examples are also provided.
Finally, Section \ref{s:Discussion}
concludes with a few remarks.

\color{black}

\section{An Overview of Inferential Models (IMs)}
\label{s:BasicAndConditionalIMs}

\replace{}{
Prior-free probabilistic inference is of utmost importance in scientific research and can be traced back to R.~A.~Fisher's inverse probability or fiducial argument \citep{fisher1930inverse, fisher1973statistical}. It can even be traced back to \cite{student1908probable}, especially when Fisher's fiducial argument is viewed as a way of constructing confidence intervals \citep{neyman1941fiducial}. Recent developments along this direction include the generalized fiducial \citep{hannig2009generalized,hannig2016generalized} and confidence distribution \citep{xie2011confidence,xie2013confidence}. The Dempster--Shafer theory of belief functions \citep{dempster1968generalization,shafer1976mathematical,shafer1990perspectives,dempster2008upper,dempster2008dempster}
can also be viewed as both an extension to consider set-valued inverse mapping and a generalization to develop methods of combining information for efficient inference; see \cite{denzux201640} for an excellent review.

The development of IMs for frequency-calibrated inference was inspired by the Dempster--Shafer theory \citep{martin2010dempster,zhang2011dempster}. In particular, the IM framework makes use of the concepts of belief and plausibility functions from the theory. On a deeper mathematical level, these functions originate from the more intuitive concepts of upper and lower probabilities, which provide an intrinsic connection with well-studied imprecise probabilities \citep[c.f.][and references therein]{gong2021judicious,martin2021imprecise,liu2021comment}. As investigated recently by \cite{martin2021imprecise}, the mathematical theory of IMs is also closely related to that of imprecise probabilities. 

It is perhaps helpful to note that frequency-calibrated inference is desirable in scientific investigations \citep{denoeux2018frequency}.  Such a desirable inference cannot be obtained in general without using upper and lower probabilities. Discrete data analysis and the famous benchmark Behrens-Fisher problem \citep[][and references therein]{martin2015marginal} provide excellent supporting examples.
Below, we provide a brief review of IMs to set the stage for later introducing the method of partial conditioning.
}

\subsection{The basic IMs}


To illustrate the basic IM framework under the current problem setting, let us first consider the simple case of $n=1$, where we only make use of a single pair of observations $(t_1, y_1)$. Thus, we have
\begin{equation} \label{eq:one-point-association}
y_1 = \vartheta_0(t_1) + u_1^\star,
\end{equation}
where $u_1^\star$ represents an unobserved realization of $U_1\sim\mathcal{N}(0,1)$.
Using the terminology of 
IMs, (\ref{eq:one-point-association}) describes an underlying sampling model which involves the unknown \emph{parameter} (function) $\vartheta\in\Theta$ and generates \emph{observed data} $X_1 = (t_1, Y_1)$ using \emph{auxiliary variable} $U_1\sim\mathcal{N}(0,1)$.
Note that the auxiliary variable is unobserved but predictable, since its distribution is fully specified.

In general, constructing an IM consists of three steps,
namely, an association (A) step, a prediction (P) step, and
a combination (C) step. These three steps
are explained below in more detail in the context of the one-point nonparametric regression problem (\ref{eq:one-point-association}):
\begin{description}
\item[A-step.]
The \emph{association step} can be achieved via some function or procedure $\mathcal{F}$ as
$$ X_1=\mathcal{F}(U_1, \vartheta), \quad (X\in\mathcal{X},\ \vartheta\in\Theta,\ U\sim\mathcal{N}(0,1)). $$
This association allows for direct reasoning with the source of uncertainty $U_1$, which is missing. If $U_1 = u_1^\star$ were observed, we would then be able to obtain the best possible inference for $\vartheta$, which is given by the set-valued ``inverse" mapping
$$ \mathcal{G}: u_1^\star \rightarrow \Theta_{X_1}(u_1^\star) = \{ \vartheta\in\Theta : \vartheta(t_1) = Y_1 - u_1^\star \}. $$
\item[P-step.]
For inference on $\vartheta$, the discussion in the A-step  suggests that we should focus our attention on accurately predicting the unobserved quantity $u_1^\star$.
To predict $u_1^\star$ with a certain desired accuracy, we utilize a \emph{predictive random set}
		$\mathcal{S}$, for example,
\begin{equation}\label{eq:PRS}
S(U_1) = \{ \tilde{u}\in\mathbb{R}: -|U_1| \le \tilde{u} \le |U_1| \},\quad 
	(U_1\sim N(0,1)).
\end{equation}

\item[C-step.]
To transfer the available information about $u^\star$ to the $\vartheta$-space, the last step is to combine the information in the association, the observed $x_1 = (t_1,y_1)$, and the predictive random set $\mathcal{S}$. For this purpose, consider the expanded set
$$ \Theta_{x_1} (\mathcal{S}) = \bigcup_{u_1\in\mathcal{S}} \Theta_{x_1} (u_1) = \{\vartheta\in\Theta: |\vartheta(t_1) - y_1| \le |U_1| \}, \quad (U_1\sim\mathcal{N}(0,1)) $$
which contains those values of $\vartheta$ that are consistent with the observed data and the sampling model for at least one candidate $u_1^\star$ value $u_1\in\mathcal{S}$.
\end{description}

The random sets obtained in the C-step are in the space of the unknown parameter, and we are ready to produce uncertainty assessment for assertions of interest.
Consider an \emph{assertion} $\mathcal{A}$ about the parameter of interest $\vartheta$. The assertion $\mathcal{A}$ corresponds to a set $A\subseteq\Theta$, and acts as a hypothesis in the context of classical statistics.
To summarize the evidence in $x$ that supports the assertion $\mathcal{A}$, we evaluate the \emph{belief function} defined by
\begin{equation}\label{eq:belief}
\mathsf{bel}_x(A; \mathcal{S}) = \mathsf{P}_\mathcal{S}\{\Theta_x(\mathcal{S}) \subseteq A\,|\,\Theta_x(\mathcal{S})\neq\varnothing\},
\end{equation}
and the \emph{plausibility function} defined by
\begin{equation}\label{eq:plausibility}
\mathsf{pl}_x(A;\mathcal{S}) = 1 - \mathsf{bel}_x(A^c, \mathcal{S}) = \mathsf{P}_\mathcal{S}\{\Theta_x(\mathcal{S})\cap A\neq\varnothing\,|\,\Theta_x(\mathcal{S}) \neq \varnothing\}.
\end{equation}
The belief function is \emph{subadditive} in the sense that if $\varnothing \neq A \subseteq \Theta$, then $\mathsf{bel}(A;\mathcal{S}) + \mathsf{bel}(A^c;\mathcal{S}) \le 1$ with equality if and only if $\Theta_x(\mathcal{S})$ is a singleton with $\mathsf{P}_\mathcal{S}$-probability 1. Therefore, it follows that $\mathsf{bel}_x (A;\mathcal{S}) \le \mathsf{pl}_x(A;\mathcal{S})$, and they are also referred to as the lower and upper probabilities, respectively.
Incidentally, we note that there are continuing interests in using 
lower and upper probabilities for statistical inference
\citep{gong2021judicious,liu2021comment}.


\begin{figure}[!htb]
  \centering
  \includegraphics[width=0.75\columnwidth]{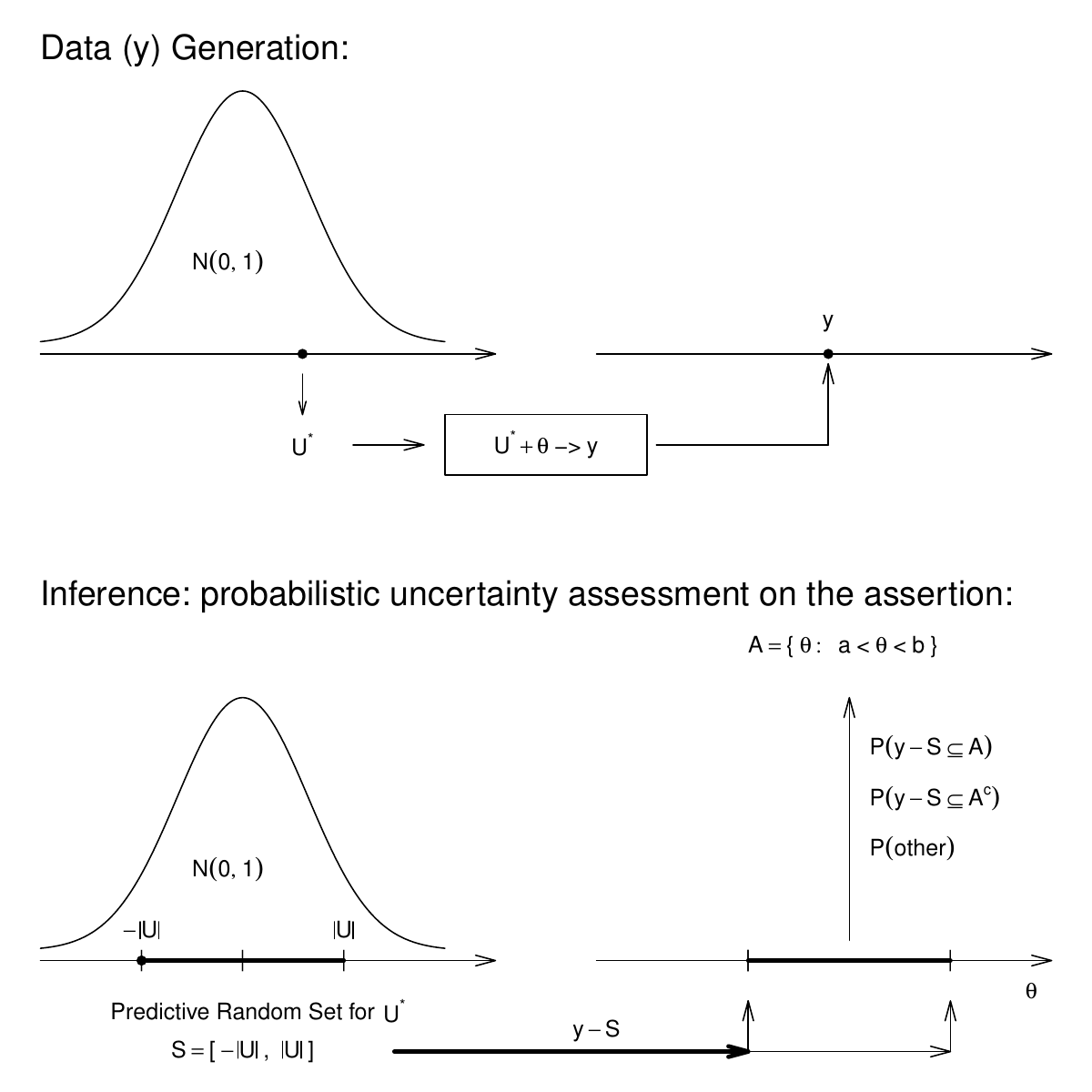}
	\caption{\replace{}{\small A pictorial illustration of IM inference about $\theta \in {\mathbb R}$ in the sampling model $Y\sim \mathcal{N}(\theta, 1)$. The A-step is given by the data generation process. The P-step uses a predictive random set $\mathcal{S}(U)=[-|U|, |U|]$, $U\sim \mathcal{N}(0,1)$, to predict the unobserved realization $U^*$ in the data generation. The C-step makes use of the random set in the $\theta$ space induced by $\mathcal{S}(U)=[-|U|, |U|]$ to compute
 $\mathsf{bel}_y(A)$, which is the probability for the truth of $A$, $P(y-S\subseteq A)$. Similarly, it computes $\mathsf{bel}_y(A^c)$, which is the probability for the falsity of $A$, $P(y-S\subseteq A^c)$ and, thereby, the plausibility $\mathsf{pl}_y(A) = 1 - P(y-S\subseteq A^c)$.
	}}
  \label{fig:basic-ims}
\end{figure}

\replace{}{
Figure \ref{fig:basic-ims} provides a pictorial illustration of IM inference about $\theta \in {\mathbb R}$ in the sampling model $Y\sim \mathcal{N}(\theta, 1)$. The A-step is given by the data generation scheme: $y = \theta+u^\star$ with $u^\star$ a realization of the random variable $U$following the known distribution $\mathcal{N}(\theta, 1)$. Inference about $\theta$, in terms of assertions 
about $\theta$, is obtained by predicting~$u^\star$ with the predictive random set ${\cal S}(U) = [-|U|, |U|]$,
$U\sim \mathcal{N}(0, 1)$.
}


To establish the \emph{frequency-calibration properties} of the IM in the context of our current problem, we briefly review several results regarding \emph{validity} as formally introduced by \cite{martin2013inferential}. 

\begin{definition}[Validity of Predictive Random Sets]\label{def:PRS}
A predictive random set $\mathcal{S}$ is $\emph{valid}$ for predicting the unobserved auxiliary variable $U$ if for each $\alpha\in(0,1)$,
\begin{equation}\label{eq:PRSvalidity}
\mathsf{P}_U \{\mathsf{P}_\mathcal{S}\{\mathcal{S}\not\owns U\}\ge1-\alpha\} \le\alpha.
\end{equation}
\end{definition}

\replace{}{
Intuitively, a valid predictive random set for a random variable $U$ should be sufficiently large so that its error rate in predicting realizations of $U$ is small enough in accordance with our familiar frequency calibration. The use of valid predictive randoms leads to frequency-calibrated inference on unknown parameters, which is defined by \cite{martin2013inferential} as follows.
}

\begin{definition}[Validity of IMs]
The IM is \emph{valid} if, for all assertions $A$ and for any $\alpha\in(0,1)$,
\begin{equation}\label{eq:valid-IM-def}
\sup_{\vartheta\in A} \mathsf{P}_{X|\vartheta}\{\mathsf{pl}_X(A;\mathcal{S})\le\alpha\} \le\alpha.
\end{equation}
\end{definition}

\replace{}{Using \eqref{eq:plausibility}, we can rewrite \eqref{eq:valid-IM-def} as
\[
\sup_{\vartheta\in A} \mathsf{P}_{X|\vartheta}\{\mathsf{bel}_X(A^c;\mathcal{S})\ge 1- \alpha\} \le\alpha.
\]
Thus, the validity of IMs guarantees that the chance of producing large degrees of belief in the truth of a false assertion is small. The particular choice of the mathematical definition is to ensure that the belief probability is frequency-calibrated. Intuitively, it essentially states that ``5\% error rate means (at most) 5\% error rate'' (e.g., the coverage-error of confidence/plausibility intervals or the Type-I error of significance testing). Incidentally, IMs are also directly related to significance testing \citep[c.f.][for details]{martin2014note}.
}

\begin{theorem}\label{thm:validity}
Suppose that the predictive random set $\mathcal{S}$ is valid, and $\Theta_x(\mathcal{S})\neq\varnothing$ with $\mathsf{P}_\mathcal{S}$-probability 1 for all $x$. Then the IM is valid.
\end{theorem}

In other words, a valid inference in the sense of
frequency calibration is obtained as long as predictive random sets
are valid to predict unobserved auxiliary variables. 
\replace{}{For more general discussion on frequency calibration of belief functions,
see \cite{denoeux2018frequency}.}

For the one-point case, we have the following result,
followed by a numerical example to illustrate
plausibility-based confidence intervals.

\begin{proposition}\label{prop:one-point-example}
The predictive random set defined by \eqref{eq:PRS} is valid and, therefore, the IM with the belief function given by \eqref{eq:belief} and the plausibility function given by \eqref{eq:plausibility} is valid.
\end{proposition}

\ifthenelse{1=1}{}{

\begin{proof}
For the predictive random set $\mathcal{S} = [-|U_1|,|U_1|],\ U_1\sim\mathcal{N}(0,1)$, we have
$$ \mathsf{P}_\mathcal{S} \{\mathcal{S}\not\owns u^\star\} = \mathsf{P} \{|U_1|<|u_1^\star|\} = 2\Phi(|u_1^\star|)-1, $$
where $\Phi(\cdot)$ denotes the standard Normal CDF.
It follows that for each $\alpha\in(0,1)$,
$$ \mathsf{P}_{U_1} \{\mathsf{P}_\mathcal{S}\{\mathcal{S}\not\owns U_1\}\ge1-\alpha\} = \mathsf{P}_{U_1} \{2\Phi(|U_1|)-1\ge 1-\alpha\} = \mathsf{P}_{U_1} \{\Phi(|U_1|)\ge 1-\alpha/2\} = \alpha, $$
which verifies that \eqref{eq:PRSvalidity} holds. Thus, we have shown that the predictive random set $\mathcal{S}$ is valid, and Theorem \ref{thm:validity} implies that the IM previously defined is valid as well.
\end{proof}

}

\begin{example}[An numerical illustration]
	\label{example:plausibity}
\rm

\begin{figure}[!htb]
  \centering
  \includegraphics[width=0.55\columnwidth]{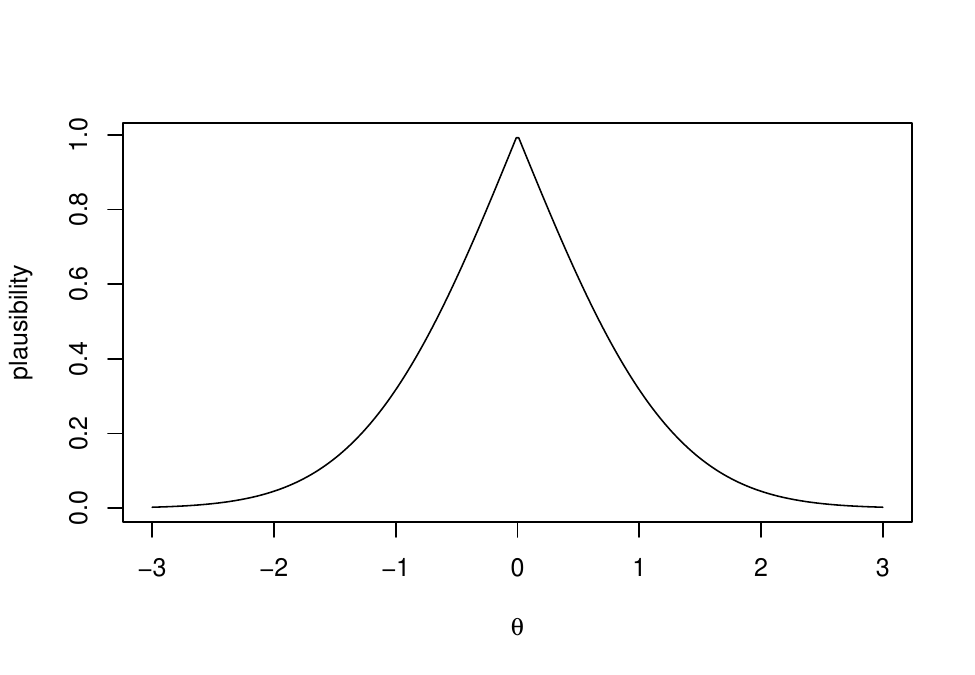}
	\caption{The plausibility function
	$\mathsf{pl}_{y_1} (A_{\theta_0};\mathcal{S})$
	defined in \eqref{eq:pl-example2.1}
	as a function of $\theta_0$ with $y_1=0$.
	}
  \label{fig:1PointPLF}
\end{figure}

	Consider the assertion $A_{\theta_0} = \{ \vartheta_0: \vartheta_0(t_1) = \theta_0 \} \subseteq \Theta_{1,1/2}$ and the predictive random set $\mathcal{S} = [-|U_1|,|U_1|],\ U_1\sim\mathcal{N}(0,1)$.
Since this assertion constitutes a singleton in the parameter space,
	the belief function $\mathsf{bel}_{y_1}(A_{\theta_0};\mathcal{S}) = 0$, but the plausibility function
\begin{equation}
	\label{eq:pl-example2.1}
 \mathsf{pl}_{y_1} (A_{\theta_0};\mathcal{S}) = \mathsf{P}_\mathcal{S} \{ \theta_0 \in \Theta_{x_1} (\mathcal{S}) \} = \mathsf{P}_{U_1} \{|U_1| \ge |\theta_0-y_1|\} = 2(1-\Phi(|\theta_{0}-y_1|)). 
\end{equation}
	Figure \ref{fig:1PointPLF} shows the graphs of $\mathsf{pl}_{y_1} (A_{\theta_0};\mathcal{S})$ as a function of $\theta_0$ for $y_1 = 0$.
\end{example}

\begin{remark}[Plausibility intervals]
	In Example \ref{example:plausibity},
	the assertion $A=\{\vartheta_0=\theta_0\},\ \vartheta_0\in\Theta$ constitutes a singleton in the parameter space. In this case, we can define the $100 (1-\alpha) \%$ \emph{plausibility region} given by
\begin{equation}\label{eq:plregion}
\Pi_x(\alpha) = \{\vartheta: \mathsf{pl}_x(\vartheta)\ge \alpha\},
\end{equation}
where $\mathsf{pl}_x(\vartheta) := \mathsf{pl}_x(\{\vartheta\};\mathcal{S})$.
It was shown in \cite{martin2013inferential} that
the plausibility region \eqref{eq:plregion} provides an exact $100 (1-\alpha) \%$ confidence interval. 
\end{remark}


\subsection{Remarks on predictive random sets}
\label{prs:efficiency}
\replace{}{It is easy to construct valid predictive random sets, for example, by simply defining a nested predictive random set
\[
{\cal S}(U )= \{u:\; b(u) \le b(U)\}, U\sim {\mathbb P}_U,
\]
with a specified boundary function $b(u)$, $u\in{\mathbb U}$, where ${\mathbb U}$ denotes the sampling space of the auxiliary variable $U$ and 
${\mathbb P}_U$, the distribution of $U$. \cite{martin2013inferential} argue for the use of nested predictive random sets. This makes sense intuitively because we want to use small set to cover unobserved $u^\star$ with large probability for the sake of efficiency.

More research is still needed to investigate the efficiency issue with respect to assertions of interest. Perhaps, investigations along the line of mathematical decision theory are a possibility. The current practice of IMs is more or less intuition based. For example, we choose centered  predictive random sets such as ${\cal S}(U) = [-|U|, |U|]$ when $U\sim\mathcal{N}(0,1))$ so that the resulting plausibility intervals are efficient in terms of interval length. Such centered predictive random sets are referred to as ``default'' predictive random sets.

}

\subsection{The marginal IMs}
\label{sec:TwoPoint-mims}
\replace{}{
As seen in the previous section for the simple one-point case, the basic IMs framework is similar to the frequentist pivotal method of constructing confidence intervals.
However, the prediction of unobserved auxiliary variables can help perform efficient inference. 
Two such methods, called
{\it conditional} IMs and {\it marginal} IMs, have been proposed in \cite{martin2015marginal,martin2015conditional} based on dimension-reduction. 
Below, we describe how marginal IMs can be used to produce
efficient inference in the special cases of the two-point problem when $B=\infty$,
where $B=M|t_1-t_2|^\gamma$.

In the $n=2$ case, we have two pairs of observations $(t_1, y_1)$ and $(t_2, y_2)$. The sampling model for this case can be written as
\begin{equation}\label{eq:2Point}
\begin{dcases}
y_1 = \vartheta_0(t_1) + u_1^\star; \\
y_2 = \vartheta_0(t_2) + u_2^\star,
\end{dcases}
\end{equation}
where $u_1^\star$ and $u_2^\star$ represent unobserved realizations of $U_1,U_2
\stackrel{iid}{\sim} \mathcal{N}(0,1)$.
The IM framework for the current problem repeats the three steps introduced in the previous section.
However, we notice that in this two-point case, we have two auxiliary variables $U_1, U_2$, which would need to be predicted using a two-dimensional predictive random set.
As \cite{martin2015marginal} 
pointed out, a more efficient inference procedure can be obtained by reducing the dimension of the auxiliary variable. 

It is easy to understand marginal IMs for the $B=\infty$ case. In this case, we have two distinct parameters $\theta_1 = \vartheta(t_1)$ and $\theta_2 = \vartheta(t_2)$. This is a problem that is discussed in \cite{martin2015marginal}. An intuitive approach is to construct a valid predictive random set, denoted by ${\mathcal S}_{U_1, U_2}$, for $(u_1^\star, u_2^\star)$ in such a way that the resulting inference on, for example, $\theta_2$ is efficient. By efficient in this case, we mean that the plausibility intervals for $\theta_2$ are as small as possible. It is easy to see that this is achieved if the projection of ${\mathcal S}(U_1, U_2)$ to the space of $u_2^\star$ is minimized, that is, ${\mathcal S}(U_1, U_2)$ has the form
\[
{\mathcal S}(U_1, U_2)= \{(u_1, u_2):\; |u_2| \leq |U_2|\}, \qquad U_1, U_2 \stackrel{iid}{\sim}\mathcal{N}(0,1).
\]
Thus, this two-dimensional predictive random set becomes effectively a one-dimensional predictive random set 
\[
{\mathcal S}(U_2) = \{(u_1, u_2):\; |u_2| \leq |U_2|\}, \qquad U_2 \sim \mathcal{N}(0,1).
\]
for $u_1^\star$. Marginal inference about $\theta_2 = \vartheta(t_2)$ proceeds with the second equation in \eqref{eq:2Point} as the association with the predictive random set ${\mathcal S}(U_2)$.
See \cite{martin2015marginal} for more discussion on marginal IMs and its more sophisticated methods.

It is seen that the above marginal inference about $\theta_2$ effectively ignores the first equation in \eqref{eq:2Point}. Although his marginal inference is still valid in the case of $B\neq \infty,$ it can be improved because the constraint $|\vartheta(t_1)-\vartheta(t_2)|\leq B$
provides a connection that allows for the observed information of $y_1$ to be used.
}
\subsection{The conditional IMs}
\label{sec:TwoPoint}

\replace{}{
Here we explain conditional IMs 
for efficient inference in the special $B=0$ case of the two-point problem.
A new method of taking the strengths of conditional IMs and marginal
IMs is proposed for the general case of the two-point problem in
Section \ref{s:PartialConditioning}.

In the $n=2$ case, we have two pairs of observations $(t_1, y_1)$ and $(t_2, y_2)$ with the sampling model given by the association \eqref{eq:2Point}. 
Again, the basic IM framework for the current problem repeats the three steps introduced in the previous section.
However, we notice that in this two-point case, we have two auxiliary variables $U_1, U_2$, which would need to be predicted using a two-dimensional predictive random set.
As \cite{martin2015conditional} 
pointed out, a more efficient inference procedure can be obtained by reducing the dimension of the auxiliary variable.
In particular, if some functions of the original auxiliary variable are fully observed, we can condition on the fully observed information to sharpen our prediction of the unobserved $u_2^\star$.
}

Notice that from \eqref{eq:2Point} we can obtain
\begin{equation}\label{eq:2PointDiff}
y_2-y_1 = \vartheta_0(t_2)-\vartheta_0(t_1) + u^\star_2-u^\star_1,
\end{equation}
which motivates us to introduce a new auxiliary variable $V=U_2-U_1\sim\mathcal{N}(0,2)$.
By the H\"older condition \eqref{eq:Holder}, we have
\begin{equation}\label{eq:2PointHolder}
|\vartheta_0(t_1)-\vartheta_0(t_2)| \le M|t_1-t_2|^\gamma.
\end{equation}
Letting $B:=M|t_1-t_2|^\gamma$, from \eqref{eq:2PointDiff} and \eqref{eq:2PointHolder} we obtain
\begin{equation}\label{eq:constraint-U12}
v^\star = u_2^\star - u_1^\star = y_2 - y_1 - (\vartheta_0(t_2)-\vartheta_0(t_1)) \in [y_2-y_1-B, y_2-y_1+B].
\end{equation}
For each fixed $v\in [y_2-y_1-B, y_2-y_1+B]$, the conditional distributions
of a linear function of $U_1$ and $U_2$, for example, $U_2$, given $V=v$ can be derived as follows:
\begin{equation}
(U_2,V)\sim\mathcal{N}_2 \left( \left[\begin{array}{c} 0 \\ 0 \end{array}\right], \left[\begin{array}{cc} 1 & 1 \\  1 &  2 \end{array}\right] \right)  \quad \Longrightarrow \quad U_2|V=v \sim \mathcal{N}\left(\frac{v}{2}, \frac{1}{2}\right).  \label{eq:2PointCondU2}
\end{equation}
Thus, these conditional distributions, sharper than the corresponding marginal distributions, motivate us to construct more efficient predictive random sets for
more efficient inference about $(\vartheta_0(t_1), \vartheta_0(t_2))$.


For the special case of $B=0$, the equation \eqref{eq:2PointDiff} simplifies to
$$ y_2-y_1 = u_2^\star-u_1^\star. $$
Most important, $u_2^\star-u_1^\star$ in this extreme case is fully
observed and, therefore, can be easily used to predict $u_1^\star$
and $u_2^\star$. Formally,
we can construct a \emph{conditional inferential model} as proposed by 
\cite{martin2015conditional}.
Similar to basic IMs, conditional IMs also have their three steps:
\begin{description}
\item[A-step.] Under the original IM framework, the association step is achieved via the \emph{baseline association}
	\eqref{eq:2Point}.
In the case of $B=0$, this baseline association can be decomposed as
\begin{align}
\begin{dcases}
y_2-y_1 = u_2^\star-u_1^\star; \\
y_2 = \vartheta(t_2)+u_2^\star.
\end{dcases}
\end{align}
This decomposition immediately suggests an alternative association. Let $\mathsf{P}_{U_2|y_2-y_1}$ denote the conditional distribution of $U_2$, given $U_2-U_1=y_2-y_1$. Since $y_2-y_1$ does not provide information on the parameter $\vartheta$, we can establish a new association
\begin{equation}\label{eq:2PointCondAssoc}
Y_2 = \vartheta(t_2)+\tilde{U}_2, \quad (\tilde{U}_2\sim\mathsf{P}_{U_2|y_2-y_1}),
\end{equation}
which is referred to as the \emph{conditional association}.
Using \eqref{eq:2PointCondAssoc}, we can associate the observed information $Y_2$ and the parameter $\vartheta$ with the new auxiliary variable $\tilde{U}_2\sim\mathsf{P}_{U_2|Y_2-Y_1}$ to get the collection of sets
$$ \Theta_{y_2}(\tilde{u}_2^\star) = \{\vartheta\in\Theta: y_2 = \vartheta(t_2) + \tilde{u}_2^\star\}. $$

\item[P-step.]
Fixing the observed value $Y_2-Y_1=y_2-y_1$, we predict the unobserved value $\tilde{u}_2^\star$ of $\tilde{U}_2$ with a \emph{conditionally valid} predictive random set $\mathcal{S}_2\sim\mathsf{P}_{\mathcal{S}_2|y_2-y_1}$.
Notice that equation \eqref{eq:2PointCondU2} implies that
		$ \tilde{U}_2\sim\mathcal{N}\left(\frac{y_2-y_1}{2},\frac{1}{2}\right), $
which gives rise to the default predictive random set (see Section \ref{prs:efficiency} given by
$$ \mathcal{S}_2 = \left\{\tilde{u}_2: \left|\tilde{u}_2-\frac{y_2-y_1}{2}\right| \le \left|\tilde{U}_2-\frac{y_2-y_1}{2}\right|\right\}. $$

\item[C-step.] We combine the results of the association and prediction steps to get
$$ \Theta_{y_2}(\mathcal{S}_2) = \bigcup_{u_2\in\mathcal{S}_2} \Theta_{y_2}(u_2) = \left\{\vartheta\in\Theta: \left|\vartheta(t_2)-\frac{y_1+y_2}{2}\right| \le \left|\tilde{U}_2-\frac{y_2-y_1}{2}\right|\right\}, $$
where $ \tilde{U}_2\sim\mathcal{N}((y_2-y_1)/2,1/2)$.
For any assertion $A\subseteq\Theta$, the corresponding \emph{conditional belief function} is given by
\begin{equation}\label{eq:cond-bel}
\mathsf{cbel}_{y_2|y_2-y_1}(A;\mathcal{S}_2) = \mathsf{P}_{\mathcal{S}_2|y_2-y_1} \{\Theta_{y_2}(\mathcal{S}_2)\subseteq A\,|\,\Theta_{y_2}(\mathcal{S}_2)\neq\varnothing\},
\end{equation}
and the \emph{conditional plausibility function} is given by
\begin{equation}\label{eq:cond-pl}
\mathsf{cpl}_{y_2|y_2-y_1}(A;\mathcal{S}_2) = 1-\mathsf{cbel}_{y_2|y_2-y_1}(A^c;\mathcal{S}_2) = \mathsf{P}_{\mathcal{S}_2|y_2-y_1}\{\Theta_{y_2}(\mathcal{S}_2)\cap A\neq\varnothing\,|\,\Theta_{y_2}(\mathcal{S}_2)\neq\varnothing\}.
\end{equation}
\end{description}

\begin{figure}[!htb]
  \centering
  \includegraphics[width=0.5\columnwidth]{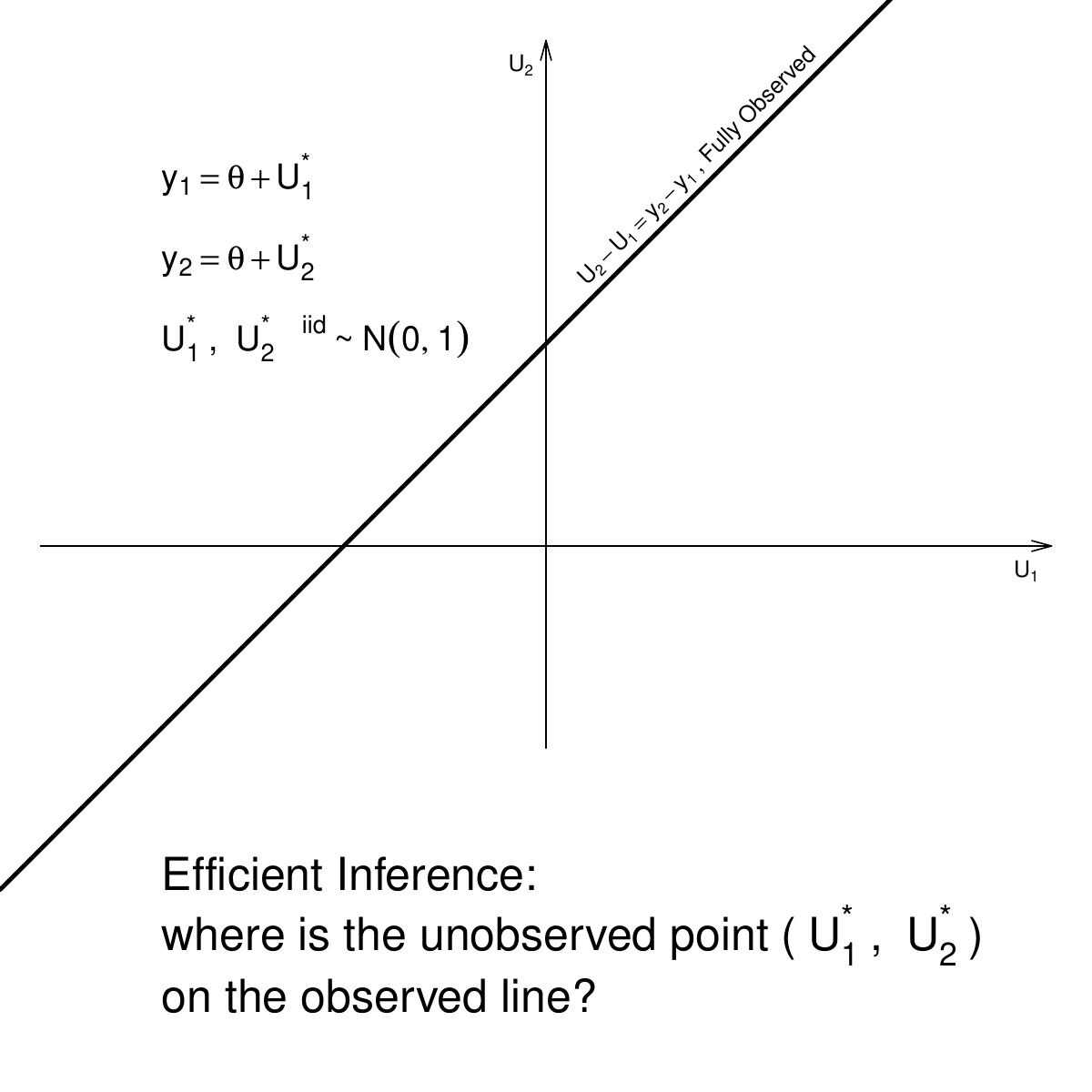}
	\caption{\replace{}{\small A pictorial illustration of conditional IM inference about $\theta \in {\mathbb R}$ in the sampling model $Y_i\stackrel{i.i.d.}{\sim} \mathcal{N}(\theta, 1)$, $i=1,2$. What is important is that $U_1-U_2=y_1-y_2$ is fully observed and, thereby, it is unnecessary to predict $(U_1^*, U_2^*)$ in the two-dimensional space of $(U_1, U_2)$. The observed information is used to further improve prediction accuracy.
	}}
  \label{fig:cond-ims}
\end{figure}

\replace{}{
Figure \ref{fig:cond-ims} provides a pictorial illustration of the above conditional IM for $\theta=\vartheta(t_1)=\vartheta(t_2)$ when $B=0$.
}
Note that, for any singleton assertion $A = \{\vartheta_0\}, \vartheta_0\in\Theta$,
the conditional belief function $\mathsf{cbel}_{y_2|y_2-y_1}(A;\mathcal{S}_2) = 0$, while the conditional plausibility function takes the form of
\[ \mathsf{cpl}_{y_2} (A;\mathcal{S}_2)  = \mathsf{P}_{\mathcal{S}_2|y_2-y_1} \{ \vartheta_0 \in \Theta_{y_2} (\mathcal{S}_2) \}
	= 2\left[1-\Phi\left(\sqrt{2}\left|\vartheta_0(t_2)-\frac{y_1+y_2}{2}\right|\right)\right].
	\]

\section{Partial Conditioning}
\label{s:PartialConditioning}
\subsection{The general $B=M|t_2-t_1|$ case and its challenges}

When $B$ is large enough, neither observation can provide information on the location
of the other observation. In this case, we may simply choose to make inference
based on associated single observation alone. The idea is supported
by the method of marginal IMs \citep{martin2015marginal}.
However, for the general case with $0 < B < \infty$,
it remains challenging to construct efficient IMs.

A seemingly attractive approach is the fiducial-type approach.
In this case, inference is obtained using the conditional distribution of
$U_2$ given $U_2 - U_1$ in the constrained interval
by taking the distribution of $U_2 - U_1$ as its original distribution
restricted to the constrained interval in \eqref{eq:constraint-U12}, {\it i.e.}, $[y_2-y_1-B, y_2-y_1+B]$.
Since this constrained interval is implicitly dependent on
the unobserved realizations of $U_2$ and $U_1$,
such a probability operation is questionable.
In fact, it generates an uncertainty assessment
without the guarantee of the desired frequency calibration.
\replace{}{More precisely,
this is because in repeated experiments for evaluating desired frequency calibration, the interval $[y_2-y_1-B, y_2-y_1+B]$ varies from experiment to experiment and depends on the experiment-specific realizations of $U_2 - U_1$. Therefore, the use of the distribution of $U_2 - U_1$ as its original distribution restricted to the interval $[y_2-y_1-B, y_2-y_1+B]$ is purely subjective and has 
no intended frequency interpretation mathematically.
 See \cite{liu2015frameworks} for more discussion of this problem for fiducial in general.
}

A conservative but simple approach to making valid inference is to
weaken the conditional belief function \eqref{eq:cond-bel} by taking its
infimum and the conditional plausibility function \eqref{eq:cond-pl} by taking
its supremum over all possible values of the conditioning variable
$U_2 - U_1$.  While valid inference is produced,
this leads to inefficient inference, which can be even worse than that based on
a single data point alone.

\subsection{A new approach based on partial conditioning}

We aim to construct nested predictive random sets for predicting $U_2$ 
based on a distribution of the form
\begin{equation}\label{eq:partial-cond}
	\mathcal{N}\left( \frac{\lambda}{2} v^\star, 1-\lambda+\frac{\lambda^2}{2} \right), 
\end{equation}
where $\lambda\in[0, 1]$ is some function of $B$ and $v^\star$ is a realization of $V=U_2-U_1$. 
\replace{}{The objective is to be able to make valid inference about $\vartheta(t_2)$ by predicting $u_2^\star$ using \eqref{eq:partial-cond} based on the fact that $U_2-U_1$ is known to lie in the constrained interval. For this reason, we call \eqref{eq:partial-cond} a predictive distribution.}
Note that the conditional and marginal IMs are obtained as two extreme cases with $\lambda=1$ for fully conditional inference and $\lambda=0$ for completely marginal inference.
The case of $0<\lambda < 1$ thus resembles an approach we refer to as
{\it partial} conditioning or, more exactly, partial
regression of, say, $U_2$ on $V=U_2 -U_1$.
The intuition and theoretical support for the use of (\ref{eq:partial-cond}) for IM-based inference are explained below\replace{}{, with a remark on its connection with the familiar method of shrinkage estimation.}

The key idea is to take into account the strengths of both conditional and marginal IMs.
That is, we consider conditional IMs but conservative for validity
when $B$ is small, and marginal IMs when conditioning-based conservative IMs
are inefficient. Technically, this amounts to utilizing the
partial regression
\[
	U_2 = \frac{\lambda}{2} V + \varepsilon
\]
to predict $U_2$ through prediction of $\varepsilon$ in such a way
that the prediction of $\varepsilon$
is valid marginally, and thereby prediction of $U_2$
is valid conditionally, given $V$ in some interval.
Marginally, $\varepsilon$ is normal with mean zero and variance
$$ \mbox{Var}(\varepsilon)= \text{Var}\left(U_2-\frac{\lambda}{2}V\right)
= \left(1-\frac{\lambda}{2}\right)^2 + \left(\frac{\lambda}{2}\right)^2 = 1-\lambda + \frac{\lambda^2}{2}. $$
This leads to the use of (\ref{eq:partial-cond}) as a valid inference.

To predict $U_2$, we utilize the predictive random sets
$$ \mathcal{S}_{v^*} = \left\{ z: \left|z-\frac{\lambda}{2} v^\star\right| \le \left|Z-\frac{\lambda}{2} v^\star\right| \right\}, \quad Z\sim\mathcal{N}\left( \frac{\lambda}{2} v^\star, 1-\lambda+\frac{\lambda^2}{2} \right) $$
when $V$ is known to be $v^*$.
For the two-point problem, $V$ is known to be some $v^*$ in
the interval $[y_2-y_1-B,\,y_2-y_1+B]$.
Thus, the use of the conservative 
prediction random set
\[
	\mathcal{S} = \bigcup_{v^*\in [y_2-y_1-B,\,y_2-y_1+B]}\mathcal{S}_{v^*}
\]
provides valid inference.
For constructing confidence intervals, for example,
the plausibility region
\[
\bigcup_{v^\star\in[y_2-y_1-B,\,y_2-y_1+B]} \left[\,\frac{\lambda}{2} v^\star - z_{1-\alpha/2} \sqrt{1-\lambda+\frac{\lambda^2}{2}},\ \frac{\lambda}{2} v^\star + z_{1-\alpha/2} \sqrt{1-\lambda+\frac{\lambda^2}{2}}\,\right]
	\]
covers $u_2^\star$ with probability at least $100(1-\alpha)\%$,
for all $\lambda \in[0,1]$.

Note that the length of the plausibility interval is given by
\begin{equation}\label{eq:partial-CI-width-02}
\left[\,\frac{\lambda}{2}(y_2-y_1+B) - \frac{\lambda}{2}(y_2-y_1-B)\,\right] + 2z_{1-\alpha/2} \sqrt{1-\lambda+\frac{\lambda^2}{2}} =  \lambda B + 2z_{1-\alpha/2} \sqrt{1-\lambda+\frac{\lambda^2}{2}}.
\end{equation}
Using the usual measure of efficiency in terms of interval length,
we minimize (\ref{eq:partial-CI-width-02})
over $\lambda \in[0,1]$.
This leads to the choice of $\lambda$:
$$ \hat{\lambda}_B =
\begin{dcases}
1 - \frac{B}{\sqrt{2z_{1-\alpha/2}^2-B^2}}, & \mbox{if } 0\le B < z_{1-\alpha/2}; \\
0, & \mbox{if } B\ge z_{1-\alpha/2}.
\end{dcases}
$$
The corresponding widths are equal to $B+\sqrt{2z_{1-\alpha/2}^2-B^2}$ and $2z_{1-\alpha/2}$, respectively.

\begin{figure}[!htb]
  \centering
  \includegraphics[width=0.5\columnwidth]{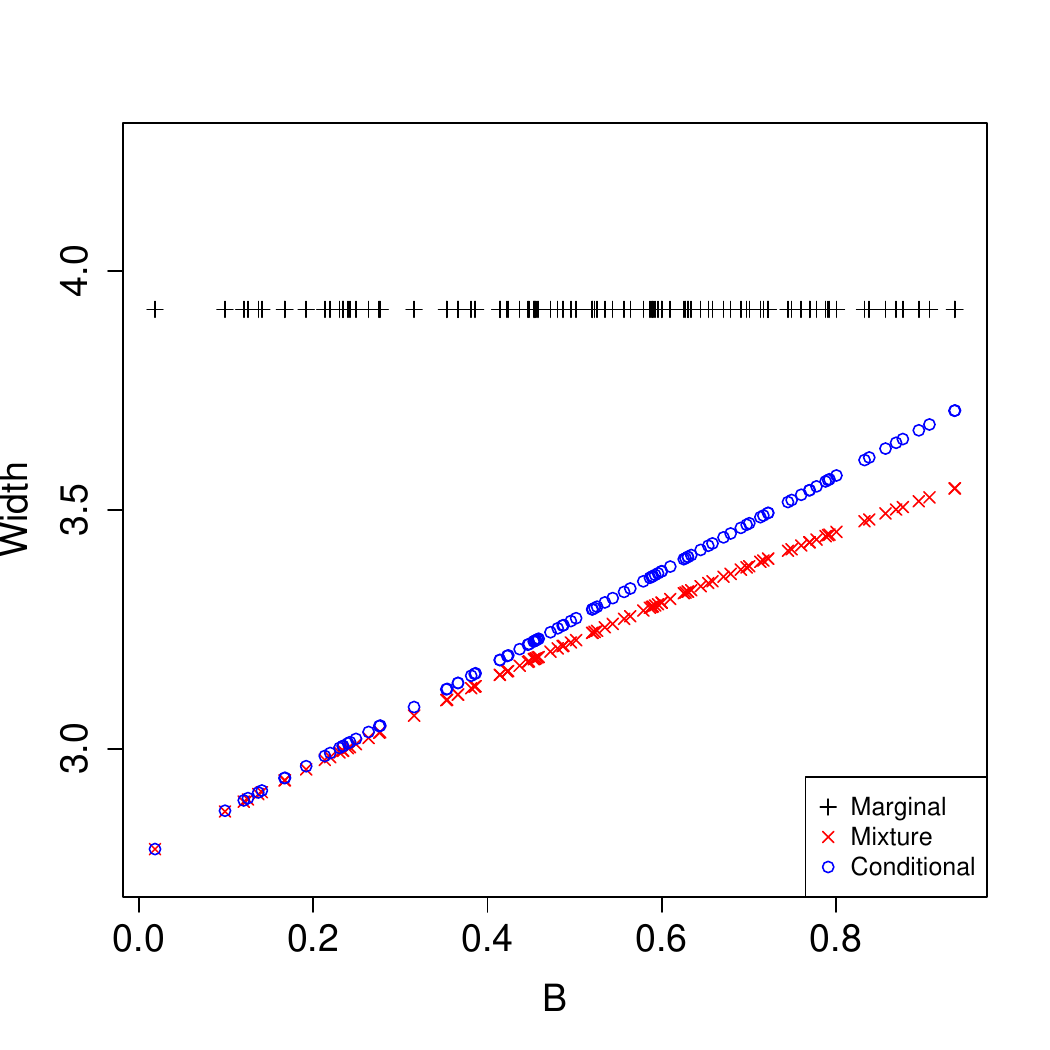} \\
  \caption{Comparison of the widths of the plausibility intervals constructed for the assertion
	{ $A(t) = \{ \vartheta_0: [0,1] \rightarrow \mathbb{R}\ |\ \vartheta_0(t) = \sqrt{t} \} \subseteq \Theta_{1,1/2}$}
	using the three different approaches:
	Marginal IMs, Partial Conditioning (``Mixture"), 
	and Conservative Conditional IMs.
	}
  \label{fig:2PMix}
\end{figure}

The optimal value of $\hat{\lambda}_B$ suggests
a full conditional IM when $B$ is small relative to
the width of the target confidence interval or the confidence level.
Figure \ref{fig:2PMix} shows the numerical results of a simulation with 100 trials. We observe that the plausibility intervals obtained
using partial conditioning are the narrowest.

\ifthenelse{1=1}{}{
\begin{lstlisting}
##### The Two-Point Case #####
set.seed(1234)

M = 1; gam = 1/2;
alpha = 0.05; z = qnorm(1-alpha/2);

N = 100 # Number of trials
W <- matrix(rep(0, N*4), N, 4)
colnames(width) <- c("Marginal", "Mixture", "Conditonal", "B")
for (i in 1:N){
  t <- sort(runif(2))
  y <- sqrt(t) + rnorm(1)
  B <- M*diff(t)^gam
  W[i,1] = 2*z
  if (B < z)
    W[i,2] = B + sqrt(2*z^2-B^2)
  else
    W[i,2] = 2*z
  W[i,3] = B + sqrt(2)*z
  W[i,4] = B
}

plot(W[,4], W[,1], ylim=c(2.75,4.25), xlab="B", ylab="Width", type='n')
points(W[,4], W[,1], pch=3, col="black")
points(W[,4], W[,2], pch=4, col="red")
points(W[,4], W[,3], pch=1, col="blue")
legend("bottomright", c("Marginal", "Mixture", "Conditional"),
       pch=c(3,4,1), col=c("black", "red", "blue"))
\end{lstlisting}
}

\subsection{Partial conditioning versus shrinkage estimation}\label{ss:skrinkage}
 \replace{}{
 The use of \eqref{eq:partial-cond} for inference about $\vartheta(t_2)$ makes a connection with the familiar method of shrinkage estimation. At a high level, both partial conditioning and the method of shrinkage estimation use information in $y_1$ about $\vartheta(t_2)$. Here we compare the two methods in terms of point estimation. The IM method of point estimation is the set of parameter values with maximum plausibility. 
 
 The basic idea of partial conditioning is to improve inference about the mean $\vartheta(t_2)$ of $y_2$ by making use of partial information of the observed $y_1$ and the constraint on the distance between $\vartheta(t_2)$ and the mean $\vartheta(t_1)$ of $y_1$. The particular construction \eqref{eq:partial-cond} is proposed for this purpose. As a result, the maximum plausibility estimate of $\vartheta(t_2)$ is the set given by
 \begin{equation}\label{eq:nonshrinkage-estimate}
 \left[\frac{(2-\lambda_B)y_2 +\lambda_B y_1}{2} - \frac{\lambda_B}{2}B,
 \frac{(2-\lambda_B)y_2 +\lambda_B y_1}{2} + \frac{\lambda_B}{2}B\right]
 \end{equation}
 when the predictive random set ${\cal S}_{v^\star}$ is used.
 
 If the inequality constraint on the distance between $\vartheta(t_2)$ and the mean $\vartheta(t_1)$ is replaced by the probabilistic
 condition $\vartheta(t_1)|\vartheta(t_2)\sim N(\vartheta(t_2), \tau^2)$
 for some $\tau>0$, that is, $y_1|\vartheta(t_2)\sim N(\vartheta(t_2), \tau^2+1)$, the maximum likelihood estimate of  $\vartheta(t_2)$ is then given by the shrinkage estimate
 \begin{equation}\label{eq:shrinkage-estimate}
 \frac{(1+\tau^2)y_2+y_1}{(1+\tau^2)+1} \,.  
 \end{equation}
 
 It is interesting to see that the shrinkage estimate \eqref{eq:shrinkage-estimate} is the center of the maximum plausibility set \eqref{eq:nonshrinkage-estimate} when $\tau^2$ is taken to be $2(1-\lambda_B)/\lambda_B$.
In general, the maximum plausibility estimate  \eqref{eq:nonshrinkage-estimate} is different from the shrinkage estimate \eqref{eq:shrinkage-estimate}, except for
the special case of $B=0$, where the maximum plausibility estimate \eqref{eq:nonshrinkage-estimate} with $\lambda_B=0$
is the shrinkage estimate derived with the corresponding assumption $\tau^2=0$.

}

\section{The general $n$ case}
\label{s:GeneralCase}

In this section, we extend the proposed partial conditioning approach for the $n=2$ case to the general $n$ case. In order to be easily comprehensible, 
we provide detailed investigation on the three-point ($n=3$) case in Section \ref{sec:3Point} with an illustrative numerical example. The general $n$ case follows as a straightforward generalization of the three-point problem and is summarized in Section \ref{sec:nPoint}, followed by numerical examples. \replace{}{Although the technical details are somewhat tedious, the main results for applications are given by the plausibility intervals \eqref{eq:plint-n}. Finally, we perform a simulation study investigating the asymptotic performance of partial conditioning.
}


\subsection{The $n=3$ case}\label{sec:3Point} 

For the case of $n=3$, we have three pairs of observations $(t_i, y_i), i=1,2,3$ 
with the data generation model:
\begin{equation}\label{eq:3Point}
\begin{dcases}
y_1 = \vartheta_0(t_1) + u_1^\star; \\
y_2 = \vartheta_0(t_2) + u_2^\star; \\
y_3 = \vartheta_0(t_3) + u_3^\star
\end{dcases}
\end{equation}
where $u_1^\star, u_2^\star, u_3^\star$ represent unobserved realizations of $U_1,U_2,U_3 \stackrel{iid}{\sim} \mathcal{N}(0,1)$.
The problem of interest is to make a (marginal) inference on
$\vartheta_0(t_i), i=1,2,3$.
Without loss of generality, assume that $t_1\le t_2\le t_3$.
To proceed, we first find the relevant conditional distributions
in Section~\ref{subsubsec: Cond}, and then
provide the optimal solution under a partial conditioning framework
in Section~\ref{subsubsec: Optim}.

\subsubsection{Conditional Distributions}\label{subsubsec: Cond}
The system of pairwise differences from \eqref{eq:3Point}, that is,
\begin{align}\label{eq:3PointDiff}
\begin{dcases}
y_2-y_1 = \vartheta_0(t_2)-\vartheta_0(t_1) + u^\star_2-u^\star_1; \\
y_3-y_1 = \vartheta_0(t_3)-\vartheta_0(t_1) + u^\star_3-u^\star_1; \\
y_3-y_2 = \vartheta_0(t_3)-\vartheta_0(t_2) + u^\star_3-u^\star_2
\end{dcases}
\end{align}
motivates us to introduce the potential conditioning auxiliary variables
$$ V_{21}=U_2-U_1,\ V_{31}=U_3-U_1,\ V_{32}=U_3-U_2. $$
Clearly, we have that
$$ (V_{21},V_{31},V_{32})\sim\mathcal{N}_3 \left( \left[\begin{array}{c} 0 \\ 0 \\ 0 \end{array}\right], \left[\begin{array}{ccc} 2 & 1 & -1 \\ 1 & 2 & 1 \\ -1 & 1 & 2 \end{array}\right] \right). $$
By the H\"older condition \eqref{eq:Holder}, we have
\begin{equation}\label{eq:3PointHolder}
\begin{dcases}
|\vartheta_0(t_2)-\vartheta_0(t_1)| \le M|t_2-t_1|^\gamma; \\
|\vartheta_0(t_3)-\vartheta_0(t_1)| \le M|t_3-t_1|^\gamma; \\
|\vartheta_0(t_3)-\vartheta_0(t_2)| \le M|t_3-t_2|^\gamma.
\end{dcases}
\end{equation}
For notational convenience, define
$$ B_{ij}:=M|t_i-t_j|^\gamma,\ 1\le i,j \le 3. $$
From \eqref{eq:3PointDiff} and \eqref{eq:3PointHolder} we obtain
the observed constraints on the pairwise differences $v_{i,}^{\star}$:
\begin{equation}\label{eq:3PointBounds}
v_{ij}^\star = u_{i}^\star - u_j^\star = y_i - y_j - (\vartheta_0(t_i)-\vartheta_0(t_j)) \in [y_i-y_j-B_{ij}, y_i-y_j+B_{ij}], \ 1\le j<i\le 3.
\end{equation}

A complete characterization of the conditional distributions of the original auxiliary variables $U_1,U_2,U_3$ given the new auxiliary variables $V_{21},V_{31},V_{32}$ is provided in Propositions \ref{prop:3PointCond1} and \ref{prop:3PointCond2}, with proofs given in the appendix.
Direct applications of these propositions yield the corresponding results with respect to the parameter $\vartheta\in\Theta$, which are summarized in Corollaries \ref{cor:3PointCond1} and \ref{cor:3PointCond2}.

\begin{proposition}\label{prop:3PointCond1}
For any fixed $v_{ij} = u_i - u_j$ with $i>j$, the conditional distributions of $U_1, U_2, U_3$ given $V_{ij} = v_{ij}$ take the following forms:
\begin{align}
U_1|v_{21} \sim \mathcal{N}\left(\frac{u_1-u_2}{2}, \frac{1}{2}\right), & \
U_1|v_{31} \sim \mathcal{N}\left(\frac{u_1-u_3}{2}, \frac{1}{2}\right), \
U_1|v_{32} \sim \mathcal{N}(0,1); \\
U_2|v_{21} \sim \mathcal{N}\left(\frac{u_2-u_1}{2}, \frac{1}{2}\right), & \
U_2|v_{32} \sim \mathcal{N}\left(\frac{u_2-u_3}{2}, \frac{1}{2}\right), \
U_2|v_{31} \sim \mathcal{N}(0,1); \\
U_3|v_{31} \sim \mathcal{N}\left(\frac{u_3-u_1}{2}, \frac{1}{2}\right), & \
U_3|v_{32} \sim \mathcal{N}\left(\frac{u_3-u_2}{2}, \frac{1}{2}\right), \
U_3|v_{21} \sim \mathcal{N}(0,1).
\end{align}
\end{proposition}

\ifthenelse{1=1}{}{

\begin{proof}
The conditional distributions for $U_1$ can be derived as follows:
\begin{align*}
\begin{dcases}
(U_1,V_{21})\sim\mathcal{N}_2 \left( \left[\begin{array}{c} 0 \\ 0 \end{array}\right], \left[\begin{array}{cc} 1 & -1 \\ -1 & 2 \end{array}\right] \right) & \quad \Rightarrow \quad U_1|V_{21}=v_{21} \sim \mathcal{N}\left(-\frac{v_{21}}{2}, \frac{1}{2}\right); \\
(U_1,V_{31})\sim\mathcal{N}_2 \left( \left[\begin{array}{c} 0 \\ 0 \end{array}\right], \left[\begin{array}{cc} 1 & -1 \\ -1 & 2 \end{array}\right] \right) & \quad \Rightarrow \quad U_1|V_{31}=v_{31} \sim \mathcal{N}\left(-\frac{v_{31}}{2}, \frac{1}{2}\right); \\
(U_1,V_{32})\sim\mathcal{N}_2 \left( \left[\begin{array}{c} 0 \\ 0 \end{array}\right], \left[\begin{array}{cc} 1 & 0 \\ 0 & 2 \end{array}\right] \right) & \quad \Rightarrow \quad U_1|V_{32}=v_{32} \sim \mathcal{N}(0,1).
\end{dcases}
\end{align*}
The conditional distributions for $U_2$ can be derived as follows:
\begin{align*}
\begin{dcases}
(U_2,V_{21})\sim\mathcal{N}_2 \left( \left[\begin{array}{c} 0 \\ 0 \end{array}\right], \left[\begin{array}{cc} 1 & 1 \\ 1 & 2 \end{array}\right] \right) & \quad \Rightarrow \quad U_2|V_{21}=v_{21} \sim \mathcal{N}\left(\frac{v_{21}}{2}, \frac{1}{2}\right); \\
(U_2,V_{31})\sim\mathcal{N}_2 \left( \left[\begin{array}{c} 0 \\ 0 \end{array}\right], \left[\begin{array}{cc} 1 & 0 \\ 0 & 2 \end{array}\right] \right) & \quad \Rightarrow \quad U_2|V_{31}=v_{31} \sim \mathcal{N}(0,1); \\
(U_2,V_{32})\sim\mathcal{N}_2 \left( \left[\begin{array}{c} 0 \\ 0 \end{array}\right], \left[\begin{array}{cc} 1 & -1 \\ -1 & 2 \end{array}\right] \right) & \quad \Rightarrow \quad U_2|V_{32}=v_{32} \sim \mathcal{N}\left(-\frac{v_{32}}{2}, \frac{1}{2}\right).
\end{dcases}
\end{align*}
The conditional distributions for $U_3$ can be derived as follows:
\begin{align*}
\begin{dcases}
(U_3,V_{21})\sim\mathcal{N}_2 \left( \left[\begin{array}{c} 0 \\ 0 \end{array}\right], \left[\begin{array}{cc} 1 & 0 \\ 0 & 2 \end{array}\right] \right) & \quad \Rightarrow \quad U_3|V_{21}=v_{21} \sim \mathcal{N}(0,1); \\
(U_3,V_{31})\sim\mathcal{N}_2 \left( \left[\begin{array}{c} 0 \\ 0 \end{array}\right], \left[\begin{array}{cc} 1 & 1 \\ 1 & 2 \end{array}\right] \right) & \quad \Rightarrow \quad U_3|V_{31}=v_{31} \sim \mathcal{N}\left(\frac{v_{31}}{2}, \frac{1}{2}\right); \\
(U_3,V_{32})\sim\mathcal{N}_2 \left( \left[\begin{array}{c} 0 \\ 0 \end{array}\right], \left[\begin{array}{cc} 1 & 1 \\ 1 & 2 \end{array}\right] \right) & \quad \Rightarrow \quad U_3|V_{32}=v_{32} \sim \mathcal{N}\left(\frac{v_{32}}{2}, \frac{1}{2}\right).
\end{dcases}
\end{align*}
The proof is complete by recalling that $v_{ij} = u_i-u_j$ for $1\le j<i\le 3$.
\end{proof}

}

\begin{corollary}\label{cor:3PointCond1}
Proposition \ref{prop:3PointCond1} yields the following predictive random sets:
\begin{align}
\vartheta(t_1) = Y_1-U_1 :\ & \mathcal{N}\left(\frac{Y_1+Y_2}{2} \pm \frac{B_{21}}{2},\, \frac{1}{2}\right) \ \bigcap \ \mathcal{N}\left(\frac{Y_1+Y_3}{2} \pm \frac{B_{31}}{2},\, \frac{1}{2}\right); \\
\vartheta(t_2) = Y_2-U_2 :\ & \mathcal{N}\left(\frac{Y_1+Y_2}{2} \pm \frac{B_{21}}{2},\, \frac{1}{2}\right) \ \bigcap \ \mathcal{N}\left(\frac{Y_2+Y_3}{2} \pm \frac{B_{32}}{2},\, \frac{1}{2}\right); \\
\vartheta(t_3) = Y_3-U_3 :\ & \mathcal{N}\left(\frac{Y_1+Y_3}{2} \pm \frac{B_{31}}{2},\, \frac{1}{2}\right) \ \bigcap \ \mathcal{N}\left(\frac{Y_2+Y_3}{2} \pm \frac{B_{32}}{2},\, \frac{1}{2}\right).
\end{align}
\replace{}{
where we abuse notation slightly and denote by $\mathcal{N}(\mu, \sigma^2)$ the ``default'' predictive random set (c.f.~Section~\ref{prs:efficiency}) given by $\{u: |u-\mu| \le |U - \mu|\}$ with $U\sim\mathcal{N}(\mu, \sigma^2)$.
}
\end{corollary}

\ifthenelse{1=1}{}{

\begin{proof}
We notice that the results in \eqref{eq:3PointBounds} imply the following:
\begin{align*}
u_1-u_2 = -v_{21} & \in [Y_1-Y_2-B_{21}, Y_1-Y_2+B_{21}], \\
u_1-u_3 = -v_{31} & \in [Y_1-Y_3-B_{31}, Y_1-Y_3+B_{31}]; \\
u_2-u_1 = v_{21} & \in [Y_2-Y_1-B_{21}, Y_2-Y_1+B_{21}], \\
u_2-u_3 = -v_{32} & \in [Y_2-Y_3-B_{32}, Y_2-Y_3+B_{32}]; \\
u_3-u_1 = v_{31} & \in [Y_3-Y_1-B_{31}, Y_3-Y_1+B_{31}], \\
u_3-u_2 = v_{32} & \in [Y_3-Y_2-B_{32}, Y_3-Y_2+B_{32}].
\end{align*}
Combining these results with Proposition \ref{prop:3PointCond1} completes the proof of the corollary.
\end{proof}

}

\begin{proposition}\label{prop:3PointCond2}
For any fixed $v_{i_1 j_1} = u_{i_1} - u_{j_1}$ and $v_{i_2 j_2} = u_{i_2} - u_{j_2}$ with $i_1>j_1$ and $i_2>j_2$, the conditional distributions of $U_1, U_2, U_3$ given $V_{i_1 j_1} = v_{i_1 j_1}$ and $V_{i_2 j_2} = v_{i_2 j_2}$ take the following forms:
\begin{align}
U_1|v_{21},v_{31},\ U_1|v_{21},v_{32},\ U_1|v_{31},v_{32} & \ \sim \ \mathcal{N}\left(\dfrac{2u_1-u_2-u_3}{3},\dfrac{1}{3}\right); \label{eq:3PointCond1} \\
U_2|v_{21},v_{31},\ U_2|v_{21},v_{32},\ U_2|v_{31},v_{32} & \ \sim \ \mathcal{N}\left(\dfrac{2u_2-u_1-u_3}{3},\dfrac{1}{3}\right); \label{eq:3PointCond2} \\
U_3|v_{21},v_{31},\ U_3|v_{21},v_{32},\ U_3|v_{31},v_{32} & \ \sim \ \mathcal{N}\left(\dfrac{2u_3-u_1-u_2}{3},\dfrac{1}{3}\right). \label{eq:3PointCond3}
\end{align}
\end{proposition}

\ifthenelse{1=1}{}{

\begin{proof}
The conditional distributions for $U_1$ can be derived as follows:
\begin{align*}
\begin{dcases}
(U_1,V_{21},V_{31}) \sim \mathcal{N}_3 \left( \left[\begin{array}{c} 0 \\ 0 \\ 0 \end{array}\right], \left[\begin{array}{ccc} 1 & -1 & -1 \\ -1 & 2 & 1 \\ -1 & 1 & 2 \end{array}\right] \right) & \quad \Rightarrow \quad U_1|v_{21},v_{31} \sim \mathcal{N}\left(-\dfrac{v_{21}+v_{31}}{3},\dfrac{1}{3}\right); \\
(U_1,V_{21},V_{32}) \sim \mathcal{N}_3 \left( \left[\begin{array}{c} 0 \\ 0 \\ 0 \end{array}\right], \left[\begin{array}{ccc} 1 & -1 & 0 \\ -1 & 2 & -1 \\ 0 & -1 & 2 \end{array}\right] \right) & \quad \Rightarrow \quad U_1|v_{21},v_{32} \sim \mathcal{N}\left(-\dfrac{2v_{21}+v_{32}}{3},\dfrac{1}{3}\right); \\
(U_1,V_{31},V_{32}) \sim \mathcal{N}_3 \left( \left[\begin{array}{c} 0 \\ 0 \\ 0 \end{array}\right], \left[\begin{array}{ccc} 1 & -1 & 0 \\ -1 & 2 & 1 \\ 0 & 1 & 2 \end{array}\right] \right) & \quad \Rightarrow \quad U_1|v_{31},v_{32} \sim \mathcal{N}\left(\dfrac{v_{32}-2v_{31}}{3},\dfrac{1}{3}\right).
\end{dcases}
\end{align*}
The conditional distributions for $U_2$ can be derived as follows:
\begin{align*}
\begin{dcases}
(U_2,V_{21},V_{31}) \sim \mathcal{N}_3 \left( \left[\begin{array}{c} 0 \\ 0 \\ 0 \end{array}\right], \left[\begin{array}{ccc} 1 & 1 & 0 \\ 1 & 2 & 1 \\ 0 & 1 & 2 \end{array}\right] \right) & \quad \Rightarrow \quad U_2|v_{21},v_{31} \sim \mathcal{N}\left(\dfrac{2v_{21}-v_{31}}{3},\dfrac{1}{3}\right); \\
(U_2,V_{21},V_{32}) \sim \mathcal{N}_3 \left( \left[\begin{array}{c} 0 \\ 0 \\ 0 \end{array}\right], \left[\begin{array}{ccc} 1 & 1 & -1 \\ 1 & 2 & -1 \\ -1 & -1 & 2 \end{array}\right] \right) & \quad \Rightarrow \quad U_2|v_{21},v_{32} \sim \mathcal{N}\left(\dfrac{v_{21}-v_{32}}{3},\dfrac{1}{3}\right); \\
(U_2,V_{31},V_{32}) \sim \mathcal{N}_3 \left( \left[\begin{array}{c} 0 \\ 0 \\ 0 \end{array}\right], \left[\begin{array}{ccc} 1 & 0 & -1 \\ 0 & 2 & 1 \\ -1 & 1 & 2 \end{array}\right] \right) & \quad \Rightarrow \quad U_2|v_{31},v_{32} \sim \mathcal{N}\left(\dfrac{v_{31}-2v_{32}}{3},\dfrac{1}{3}\right).
\end{dcases}
\end{align*}
The conditional distributions for $U_3$ can be derived as follows:
\begin{align*}
\begin{dcases}
(U_3,V_{21},V_{31}) \sim \mathcal{N}_3 \left( \left[\begin{array}{c} 0 \\ 0 \\ 0 \end{array}\right], \left[\begin{array}{ccc} 1 & 0 & 1 \\ 0 & 2 & 1 \\ 1 & 1 & 2 \end{array}\right] \right) & \quad \Rightarrow \quad U_3|v_{21},v_{31} \sim \mathcal{N}\left(\dfrac{2v_{31}-v_{21}}{3},\dfrac{1}{3}\right); \\
(U_3,V_{21},V_{32}) \sim \mathcal{N}_3 \left( \left[\begin{array}{c} 0 \\ 0 \\ 0 \end{array}\right], \left[\begin{array}{ccc} 1 & 0 & 1 \\ 0 & 2 & -1 \\ 1 & -1 & 2 \end{array}\right] \right) & \quad \Rightarrow \quad U_3|v_{21},v_{32} \sim \mathcal{N}\left(\dfrac{v_{21}+2v_{32}}{3},\dfrac{1}{3}\right); \\
(U_3,V_{31},V_{32}) \sim \mathcal{N}_3 \left( \left[\begin{array}{c} 0 \\ 0 \\ 0 \end{array}\right], \left[\begin{array}{ccc} 1 & 1 & 1 \\ 1 & 2 & 1 \\ 1 & 1 & 2 \end{array}\right] \right) & \quad \Rightarrow \quad U_3|v_{31},v_{32} \sim \mathcal{N}\left(\dfrac{v_{31}+v_{32}}{3},\dfrac{1}{3}\right).
\end{dcases}
\end{align*}
It can be easily verified that the three conditional distributions in each group are actually equivalent, and that the results \eqref{eq:3PointCond1} -- \eqref{eq:3PointCond3} hold.
\end{proof}

}

\begin{corollary}\label{cor:3PointCond2}
Proposition \ref{prop:3PointCond2} yields the following predictive random sets:
\begin{align*}
\vartheta(t_1) = Y_1-U_1 :\ & \mathcal{N}\left(\frac{Y_1+Y_2+Y_3}{3} \pm \frac{B_{21}+B_{31}}{3},\, \frac{1}{3}\right); \\
\vartheta(t_2) = Y_2-U_2 :\ & \mathcal{N}\left(\frac{Y_1+Y_2+Y_3}{3} \pm \frac{B_{21}+B_{32}}{3},\, \frac{1}{3}\right); \\
\vartheta(t_3) = Y_3-U_3 :\ & \mathcal{N}\left(\frac{Y_1+Y_2+Y_3}{3} \pm \frac{B_{31}+B_{32}}{3},\, \frac{1}{3}\right).
\end{align*}
\replace{}{
where we again abuse notation and denote by $\mathcal{N}(\mu, \sigma^2)$ the ``default'' predictive random set $\{u: |u-\mu| \le |U - \mu|\}$ with $U\sim\mathcal{N}(\mu, \sigma^2)$.
}
More compactly, the predictive random set for $\vartheta(t_i)$ can be written as
$$ \mathcal{N}\left(\bar{Y} \pm \bar{B}_i,\, \frac{1}{3}\right), \quad i=1,2,3. $$
where $\bar{Y}=(Y_1+Y_2+Y_3)/3, \bar{B}_i=(B_{i1}+B_{i2}+B_{i3})/3$, and we define $B_{ii}=0$ for $i=1,2,3$.
\end{corollary}

\ifthenelse{1=1}{}{

\begin{proof}
We notice that the results in \eqref{eq:3PointBounds} implies the following:
\begin{align*}
2u_1-u_2-u_3 = -(v_{21}+v_{31}) & \in [2Y_1-Y_2-Y_3-(B_{21}+B_{31}), 2Y_1-Y_2-Y_3+(B_{21}+B_{31})]; \\
2u_2-u_1-u_3 = v_{21}-v_{32} & \in [2Y_2-Y_1-Y_3-(B_{21}+B_{32}), 2Y_2-Y_1-Y_3+(B_{21}+B_{32})]; \\
2u_3-u_1-u_2 = v_{31}+v_{32} & \in [2Y_3-Y_1-Y_2-(B_{31}+B_{32}), 2Y_1-Y_2-Y_3+(B_{31}+B_{32})].
\end{align*}
Combining these results with Proposition \ref{prop:3PointCond2} completes the proof of the corollary.
\end{proof}

}

\begin{remark}
For the sake of completeness, we remark that the covariance matrices in
\begin{align*}
(U_1,V_{21},V_{31},V_{32}) \sim \left( \left[\begin{array}{c} 0 \\ 0 \\ 0 \\ 0 \end{array}\right], \left[\begin{array}{cccc} 1 & -1 & -1 & 0 \\ -1 & 2 & 1 & -1 \\ -1 & 1 & 2 & 1 \\ 0 & -1 & 1 & 2 \end{array}\right] \right); \\
(U_2,V_{21},V_{31},V_{32}) \sim \left( \left[\begin{array}{c} 0 \\ 0 \\ 0 \\ 0 \end{array}\right], \left[\begin{array}{cccc} 1 & 1 & 0 & -1 \\ 1 & 2 & 1 & -1 \\ 0 & 1 & 2 & 1 \\ -1 & -1 & 1 & 2 \end{array}\right] \right) ; \\
(U_3,V_{21},V_{31},V_{32}) \sim \left( \left[\begin{array}{c} 0 \\ 0 \\ 0 \\ 0 \end{array}\right], \left[\begin{array}{cccc} 1 & 0 & 1 & 1 \\ 0 & 2 & 1 & -1 \\ 1 & 1 & 2 & 1 \\ 1 & -1 & 1 & 2 \end{array}\right] \right),
\end{align*}
are all degenerate, due to the fact that $V_{21},V_{31},V_{32}$ are linearly dependent.
\end{remark}

\subsubsection{Prediction via partial conditioning}
\label{subsubsec: Optim}

Recall that our goal is to predict $u_2^\star$ to make a valid probabilistic inference.
Without loss of generality, assume that $t_2-t_1 \le t_3-t_2$.
Motivated by the idea of partial conditioning for the two-point problem,
we construct nested predictive random sets by introducing the 
partial regression of $u_2^\star$ on pairwise differences:
\begin{equation}\label{eq:3PMix}
\mathcal{N}\left( \lambda_1\left(\frac{u_2^\star-u_1^\star}{2}\right) + \lambda_2\left(\frac{2u_2^\star-u_1^\star-u_3^\star}{3}\right), \ \left(\frac{\lambda_1}{2}+\frac{\lambda_2}{3}\right)^2 + \left(1-\frac{\lambda_1}{2}-\frac{2}{3}\lambda_2\right)^2 + \left(\frac{\lambda_2}{3}\right)^2 \right),
\end{equation}
where the mixing proportions $\lambda_1,\lambda_2$ are functions of $B_{ij}$ satisfying $\lambda_1,\lambda_2\ge0$ and $\lambda_1+\lambda_2\le1$. Note that the variance expression arises naturally from the fact that
\begin{align*}
\text{Var}\left[U_2-\frac{\lambda_1}{2}\left(U_2-U_1\right)-\frac{\lambda_2}{3}\left(2U_2-U_1-U_3\right)\right]
= \left(\frac{\lambda_1}{2}+\frac{\lambda_2}{3}\right)^2 + \left(1-\frac{\lambda_1}{2}-\frac{2}{3}\lambda_2\right)^2 + \left(\frac{\lambda_2}{3}\right)^2.
\end{align*}

To predict $U_2$, we take the predictive random sets
$$ \mathcal{S}_{v_{21}^\star, v_{32}^\star} = \left\{ z: \left|z-\left[\lambda_1\left(\frac{v_{21}^\star}{2}\right)+\lambda_2\left(\frac{v_{21}^\star-v_{32}^\star}{3}\right)\right] \right| \le \left|Z-\left[\lambda_1\left(\frac{v_{21}^\star}{2}\right)+\lambda_2\left(\frac{v_{21}^\star-v_{32}^\star}{3}\right)\right] \right| \right\}, $$ 
where we recall that
$$ v_{ij}^\star = u_{i}^\star - u_j^\star = y_i - y_j - (\vartheta_0(t_i)-\vartheta_0(t_j)) \in [y_i-y_j-B_{ij}, y_i-y_j+B_{ij}], \ 1\le j<i\le 3, $$
and $Z$ follows the Gaussian distribution \eqref{eq:3PMix}.
The predictive random sets are clearly marginally valid, and thus the corresponding IMs are  valid. 
Therefore,
the we utilize the conservative predictive random set
\[
\mathcal{S} =
\bigcup_{\substack{v_{21}^\star\in[y_2-y_1\pm B_{21}] \\ v_{32}^\star\in[y_3-y_2\pm B_{32}]}
	}
\mathcal{S}_{v_{21}^\star, v_{32}^\star}
	\]
For constructing confidence intervals,
this suggests the plausibility region
$$ \bigcup_{\substack{v_{21}^\star\in[y_2-y_1\pm B_{21}] \\ v_{32}^\star\in[y_3-y_2\pm B_{32}]}}
\left\{ \left[\left(\frac{\lambda_1}{2}+\frac{\lambda_2}{3}\right) v_{21}^\star - \left(\frac{\lambda_2}{3}\right) v_{32}^\star\right]
\pm z_{1-\alpha/2} \sqrt{\left(\frac{\lambda_1}{2}+\frac{\lambda_2}{3}\right)^2 + \left(1-\frac{\lambda_1}{2}-\frac{2}{3}\lambda_2\right)^2 + \left(\frac{\lambda_2}{3}\right)^2} \right\}, $$
which covers $u_2^\star$ with probability at least $100(1-\alpha)\%$.
It can be easily verified that the width of this plausibility region is given by
$$ 2\left[ \left(\frac{\lambda_1}{2}+\frac{\lambda_2}{3}\right) B_{21} + \left(\frac{\lambda_2}{3}\right) B_{32} \right]
+ 2z_{1-\alpha/2} \sqrt{\left(\frac{\lambda_1}{2}+\frac{\lambda_2}{3}\right)^2 + \left(1-\frac{\lambda_1}{2}-\frac{2}{3}\lambda_2\right)^2 + \left(\frac{\lambda_2}{3}\right)^2}. $$

\begin{figure}[!htb]
  \centering
  \subfigure[Width vs $B_{12}$]{
  \includegraphics[width=0.48\columnwidth]{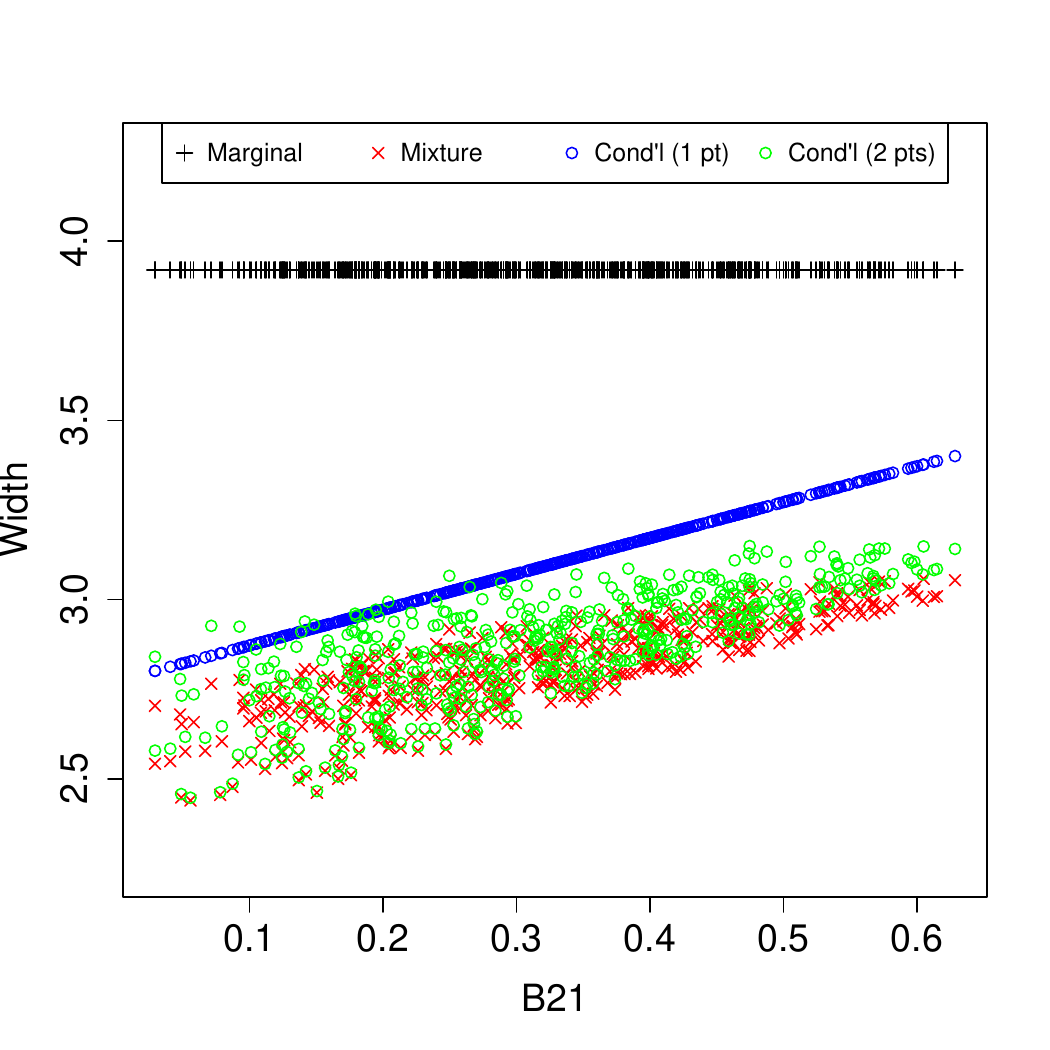}
  \label{fig:3PMix1}}
  \subfigure[Width vs $B_{23}$]{
  \includegraphics[width=0.48\columnwidth]{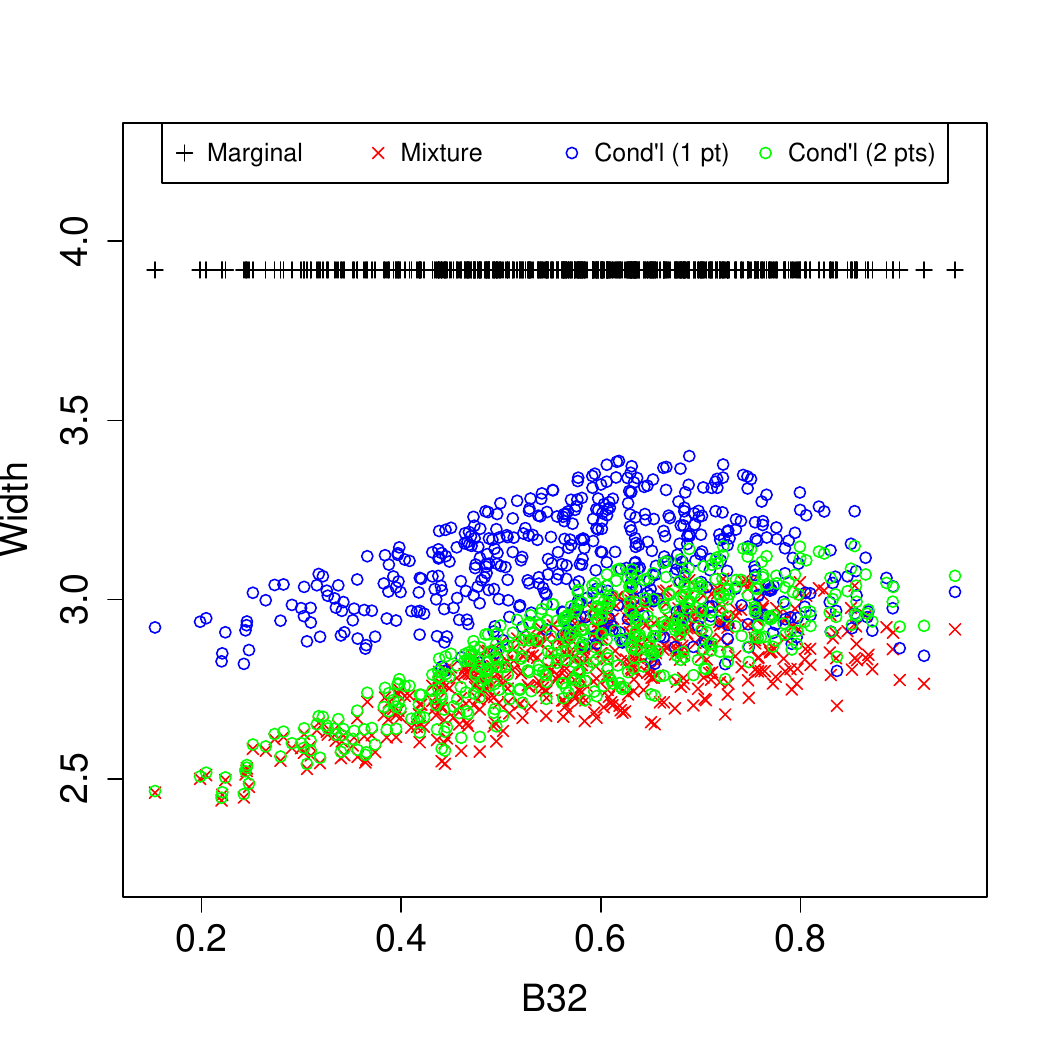}}
  \caption{Comparison of the widths of the plausibility intervals constructed for the assertion $A(t) = \{ \vartheta_0: [0,1] \rightarrow \mathbb{R}\ |\ \vartheta_0(t) = \sqrt{t} \} \subseteq \Theta_{1,1/2}$ using four different approaches: Marginal IMs, Partial Conditioning (``Mixture"), Conditional IMs on the nearest point (labeled as ``Cond'l (1 pt)" in the plots), and Conditional IMs on both points (labeled as ``Cond'l (2 pts)" in the plots).}
  \label{fig:3PMix}
\end{figure}

The optimal mixing proportions $\lambda_1,\lambda_2$ are the solutions to the constrained optimization problem
\begin{align}\label{eq:3PMixOpt}
\nonumber \min_{\lambda_1,\,\lambda_2}\quad & \left[ \left(\frac{\lambda_1}{2}+\frac{\lambda_2}{3}\right) B_{21} + \left(\frac{\lambda_2}{3}\right) B_{32} \right]
+ z_{1-\alpha/2} \sqrt{\left(\frac{\lambda_1}{2}+\frac{\lambda_2}{3}\right)^2 + \left(1-\frac{\lambda_1}{2}-\frac{2}{3}\lambda_2\right)^2 + \left(\frac{\lambda_2}{3}\right)^2} \\
\text{s.t.}\quad & \lambda_1,\lambda_2\ge0;\ \lambda_1+\lambda_2\le1.
\end{align}
Solving optimization problem \eqref{eq:3PMixOpt} analytically requires finding the roots of a set of quadratic equations, a process which can become quite burdensome.
Instead, we adopt numerical optimization methods such as the Broyden--Fletcher--Goldfarb--Shanno~(BFGS) algorithm~\citep{nocedal2006,liu2008statistical} to solve problem \eqref{eq:3PMixOpt} in order to obtain the optimal mixing proportions which minimize the length of the plausibility interval.
Figure \ref{fig:3PMix} shows the numerical results of a simulation with 500 trials.  Again, we find that the plausibility intervals obtained using optimal mixture distributions are the narrowest. 

\ifthenelse{1=1}{}{
\begin{lstlisting}
##### The Three-Point Case #####
set.seed(1234)

M = 1; gam = 1/2;
alpha = 0.05; z = qnorm(1-alpha/2);

N = 500 # Number of trials
W <- matrix(rep(0, N*6), N, 6)
colnames(W) <- c("Marginal", "Mixture", "Cond1", "Cond2", "B12", "B23")
for (i in 1:N)%{
  t <- sort(runif(3))
  y <- sqrt(t) + rnorm(3)
  B <- M*diff(t)^gam
  B <- sort(B) # Set B21 < B32
  W[i,1] = 2*z # Marginal
  # Calculate mixing proportions
  L1 <- matrix(c(1/2,1/3,0,1/3), 2, 2, byrow=T)
  L2 <- matrix(c(1/2,1/2,1/2,2/3), 2, 2, byrow=T)
  mixture.width <- function (lambda) {
    L1 <- matrix(c(1/2,1/3,0,1/3), 2, 2, byrow=T)
    L2 <- matrix(c(1/2,1/2,1/2,2/3), 2, 2, byrow=T)
    as.numeric(B%*%L1%*%lambda +
                 z*sqrt(lambda%*%L2%*%lambda+c(-1,-4/3)%*%lambda+1))
  }
  mixture.width.grad <- function (lambda) {
    t(B%*%L1) + as.numeric(z/2/sqrt(lambda%*%L2%*%lambda+c(-1,-4/3)%*%lambda+1))*((L2+t(L2))%*%lambda+c(-1,-4/3))
  }
  mixture.sol <- constrOptim(c(0.1,0.1), f=mixture.width, grad=mixture.width.grad,
                             ui=matrix(c(1,0,0,1,-1,-1),3,2,byrow=T), ci=c(0,0,-1))
  W[i,2] = 2*mixture.sol$value
  W[i,3] = B[1] + 2*z/sqrt(2) # Conditional on nearer point
  W[i,4] = 2*sum(B)/3 + 2*z/sqrt(3) # Conditional on both points
  W[i,5:6] <- B
%}

plot(W[,5], W[,1], ylim=c(2.25,4.25), xlab="B21", ylab="Width", type='n')
points(W[,5], W[,1], pch=3, col="black")
points(W[,5], W[,2], pch=4, col="red")
points(W[,5], W[,3], pch=1, col="blue")
points(W[,5], W[,4], pch=1, col="green")
legend("top", c("Marginal", "Mixture", "Cond'l (1 pt)", "Cond'l (2 pts)"),
       pch=c(3,4,1,1), col=c("black", "red", "blue", "green"), horiz=TRUE)

plot(W[,6], W[,1], ylim=c(2.25,4.25), xlab="B32", ylab="Width", type='n')
points(W[,6], W[,1], pch=3, col="black")
points(W[,6], W[,2], pch=4, col="red")
points(W[,6], W[,3], pch=1, col="blue")
points(W[,6], W[,4], pch=1, col="green")
legend("top", c("Marginal", "Mixture", "Cond'l (1 pt)", "Cond'l (2 pts)"),
       pch=c(3,4,1,1), col=c("black", "red", "blue", "green"), horiz=TRUE)
\end{lstlisting}
}

\subsection{The general $n$ case}\label{sec:nPoint}

In this section, we extend our previous discussions to the general case of $n$ observations, and provide a generic method to construct valid and efficient pointwise plausibility intervals.
Denote our sequence of observations by $(t_i, y_i),\,i=1,\cdots,n$, where it is not necessary to assume that the $t_i$'s are sorted.
With the assumption $\sigma=1$ made in Section \ref{s:Introduction},
the association \eqref{eq:model} can be written as
\begin{equation}\label{eq:nPoint}
y_i = \vartheta_0(t_i) + u_i^\star,\,i=1,\cdots,n.
\end{equation}
In the rest of this section, we outline the process of conducting inference on $U_i$ for an arbitrary $1\le i\le n$.
As a preprocessing step, we sort the $n$ observations in ascending order of their distance from the $i$-th observation,
such that $0=|t_i-t^{(i)}_0|\le |t_i-t^{(i)}_1|\le \cdots\le |t_i-t^{(i)}_{n-1}|$, correspondingly.
Denote by $ U^{(i)}_0, U^{(i)}_1, \cdots, U^{(i)}_{n-1}$
the corresponding $U$ variables, and
by $ y^{(i)}_0, y^{(i)}_1, \cdots, y^{(i)}_{n-1}$
the corresponding $y$ observations.
Notice that $U^{(i)}_0 = U_i$ holds trivially.

To predict $U_i$, we can utilize the marginal distribution $\mathcal{N}(0,1)$, and the conditional distributions $U_i|u_i-u^{(i)}_1;\ U_i|u_i-u^{(i)}_1,u_i-u^{(i)}_2;\ \cdots;\ U_i|u_i-u^{(i)}_1,u_i-u^{(i)}_2,\cdots,u_i-u^{(i)}_{n-1}$.
This motivates us to consider the partial regression of $U_i$ on the pairwise differences of the $U$ variables:
\begin{equation}\label{eq:nPMix}
\mathcal{N}\left( \sum_{k=1}^{n-1}\frac{\lambda_k}{k+1}\left[k\,u_i-\sum_{j=1}^k u^{(i)}_j\right],\,\left(1-\sum_{k=1}^{n-1}\frac{\lambda_k\,k}{k+1}\right)^2 + \sum_{j=1}^{n-1}\left(\sum_{k=j}^{n-1}\frac{\lambda_k}{k+1}\right)^2 \right),
\end{equation}
where the mixing proportions $\lambda_1,\lambda_2,\cdots,\lambda_{n-1}$ are functions of $B_{ij}$ satisfying $\lambda_1,\lambda_2,\cdots,\lambda_{n-1}\ge0$ and $\sum_{k=1}^{n-1}\lambda_k\le1$. 
The variance expression in Equation \eqref{eq:nPMix} is obtained by routine algebraic operations:
\begin{align}\label{eq:nPVar}
\nonumber & \text{Var}\left\{ U_i - \sum_{k=1}^{n-1}\frac{\lambda_k}{k+1}\left[k\,U_i-\sum_{j=1}^k U^{(i)}_j\right] \right\} \\
= & \left(1-\sum_{k=1}^{n-1}\frac{\lambda_k\,k}{k+1}\right)^2 + \text{Var}\left\{\sum_{k=1}^{n-1}\left[\frac{\lambda_k}{k+1}\sum_{j=1}^k U^{(i)}_j\right] \right\},
\end{align}
To deal with the second term on the right-hand-side of the last equation, we notice that the summation signs can be exchanged in the following manner:
$$ \sum_{k=1}^{n-1}\left[\frac{\lambda_k}{k+1}\sum_{j=1}^k U^{(i)}_j\right] = \sum_{k=1}^{n-1}\sum_{j=1}^{k} \frac{\lambda_k}{k+1}U^{(i)}_j = \sum_{j=1}^{n-1}\sum_{k=1}^{n-1} \frac{\lambda_k}{k+1}U^{(i)}_j = \sum_{j=1}^{n-1} \left(\sum_{k=j}^{n-1} \frac{\lambda_k}{k+1}\right) U^{(i)}_j. $$
Therefore, we have
$$ \text{Var}\left\{\sum_{k=1}^{n-1}\left[\frac{\lambda_k}{k+1}\sum_{j=1}^k U^{(i)}_j\right] \right\} = \text{Var}\left\{\sum_{j=1}^{n-1}\left(\sum_{k=j}^{n-1} \frac{\lambda_k}{k+1}\right) U^{(i)}_j\right\} = \sum_{j=1}^{n-1} \left(\sum_{k=j}^{n-1}\frac{\lambda_k}{k+1}\right)^2, $$
where the last step follows from the independence of $U^{(i)}_j,\,j=1,\cdots,n$.
Substituting this result into equation \eqref{eq:nPVar} yields the variance expression shown in equation \eqref{eq:nPMix}.

\begin{figure}[!htb]
  \centering
	\subfigure[$n=20,\,\sigma=1$]{
    \includegraphics[width=0.48\columnwidth]{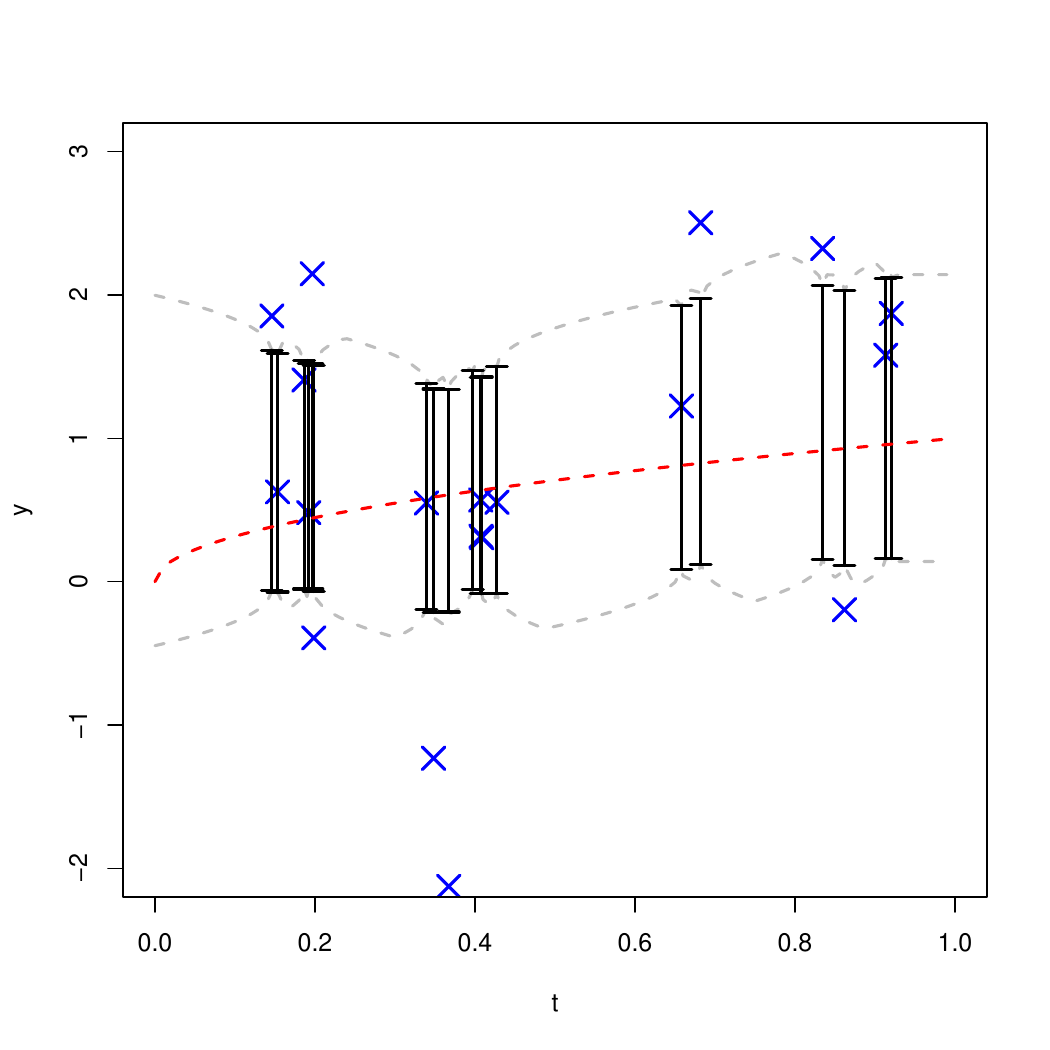}}
  \subfigure[$n=20,\,\sigma=1/2$]{
    \includegraphics[width=0.48\columnwidth]{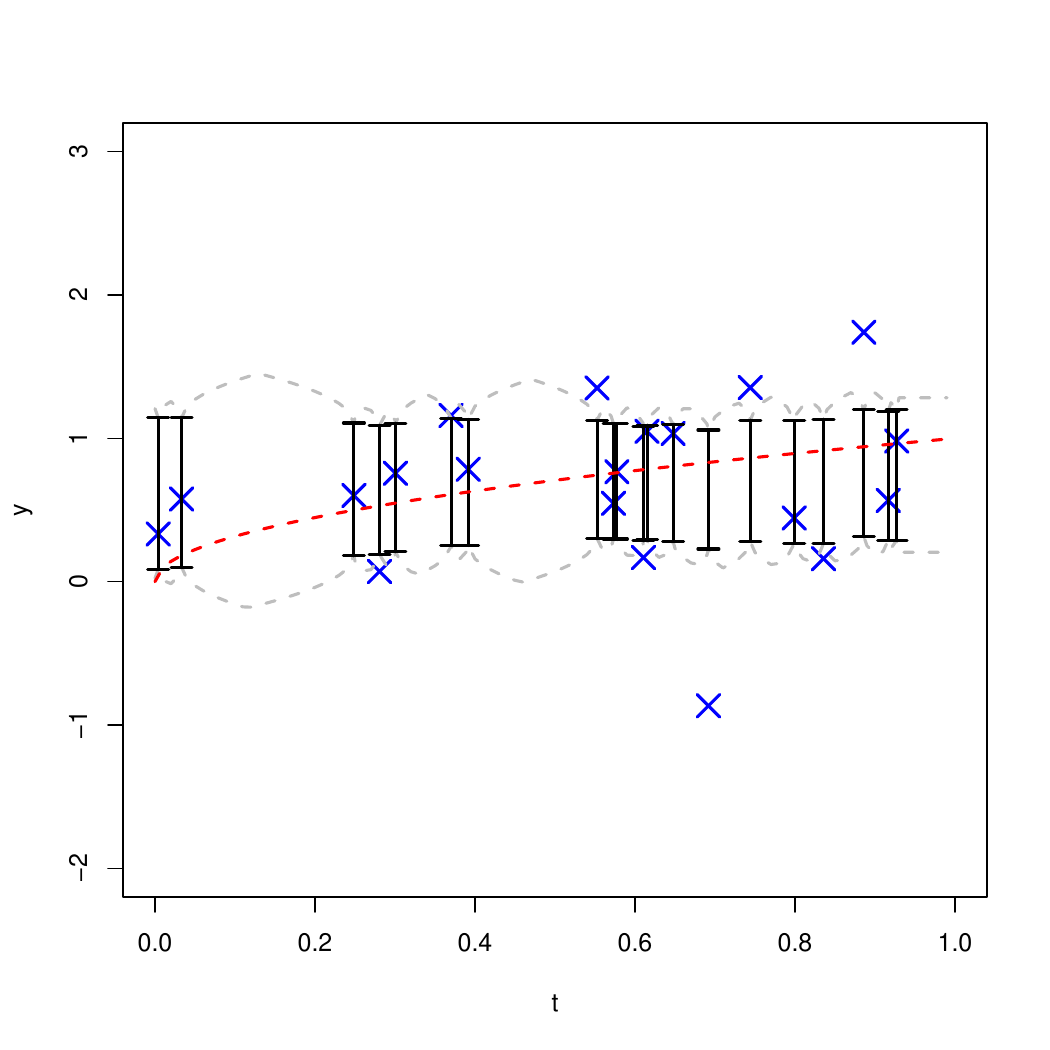}}
  \caption{Visualization of the plausibility intervals constructed using the mixture distribution for the assertion $A(t) = \{ \vartheta_0: [0,1] \rightarrow \mathbb{R}\ |\ \vartheta_0(t) = \sqrt{t} \} \subseteq \Theta_{1,1/2}$ under various settings.}
	\label{fig:n-01}
\end{figure}

\begin{figure}[!htb]
  \centering
  \subfigure[$n=15,\,M=1/2,\,\gamma=1/5$]{
    \includegraphics[width=0.48\columnwidth]{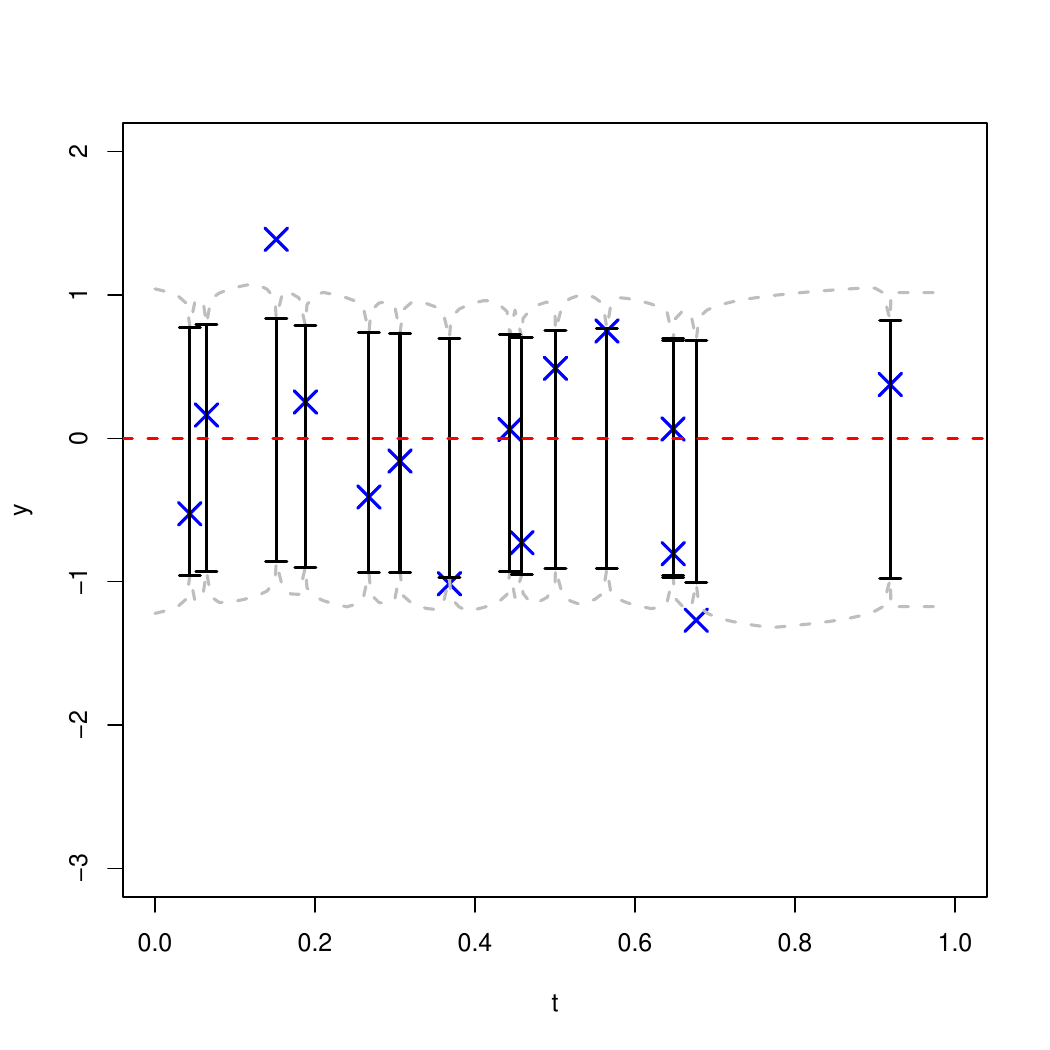}}
  \subfigure[$n=15,\,M=1/2,\,\gamma=1$]{
    \includegraphics[width=0.48\columnwidth]{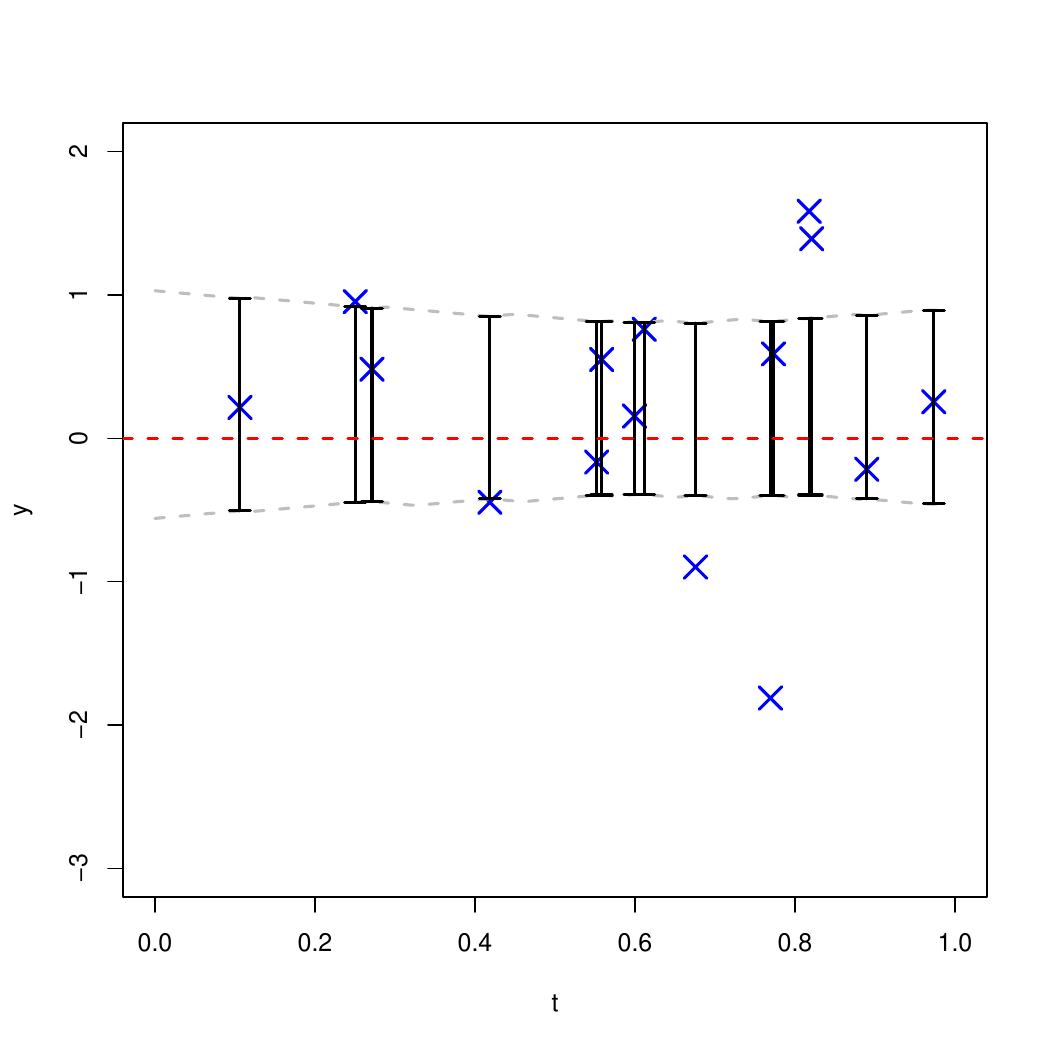}}
  \caption{Visualization of the plausibility intervals constructed using the mixture distribution for the assertion $A(t) = \{ \vartheta_0: [0,1] \rightarrow \mathbb{R}\ |\ \vartheta_0(t) = 0 \} \subseteq \Theta_{1,1/2}$ under various settings.}
	\label{fig:n-02}
\end{figure}

To predict $U_i$, we take the predictive random sets
$$ \mathcal{S}_{\{v_j^{(i)}: j=1,...,k\}} = \left\{ z: \left|z-\sum_{k=1}^{n-1}\frac{\lambda_k}{k+1} \sum_{j=1}^k v^{(i)}_j \right| \le \left|Z-\sum_{k=1}^{n-1}\frac{\lambda_k}{k+1} \sum_{j=1}^k v^{(i)}_j \right| \right\},\ Z\sim\eqref{eq:nPMix}, $$
where
$$ v^{(i)}_j := u_i-u^{(i)}_j = y_i - y^{(i)}_j - (\vartheta_0(t_i)-\vartheta_0(t^{(i)}_j)) \in [y_i-y^{(i)}_j-B^{(i)}_j, y_i-y^{(i)}_j+B^{(i)}_j], $$
with $B^{(i)}_j := |t^{(i)}_j-t_i|, \ 1\le i,j\le n$.
For valid inference, we take
the (conservative) predictive random sets
\[ \mathcal{S} = \bigcup_{\{v^{(i)}_j \in \left[y_i-y^{(i)}_j\pm B^{(i)}_j\right],j=1,...,k\}}
 \mathcal{S}_{\{v_j^{(i)}: j=1,...,k\}}.
	\]
This predictive random set is clearly marginally valid. 
Let 
\[
	\Delta_\lambda = \sqrt{\left(1-\sum_{k=1}^{n-1}\frac{\lambda_k\,k}{k+1}\right)^2 + \sum_{j=1}^{n-1}\left(\sum_{k=j}^{n-1}\frac{\lambda_k}{k+1}\right)^2},
	\]
the standard deviation of the partial regression. 
The plausibility region
$$ \bigcup_{v^{(i)}_j \in \left[y_i-y^{(i)}_j\pm B^{(i)}_j\right]} \left\{ \sum_{k=1}^{n-1}\frac{\lambda_k}{k+1} \sum_{j=1}^k v^{(i)}_j \pm z_{1-\alpha/2} 
\Delta_\lambda\right\}, $$
or equivalently,
\begin{equation}\label{eq:plint-n}
\left[ \sum_{k=1}^{n-1}\frac{\lambda_k}{k+1} \sum_{j=1}^k \left[y_i-y^{(i)}_j-B^{(i)}_j\right] - z_{1-\alpha/2} \Delta_\lambda,\ \sum_{k=1}^{n-1}\frac{\lambda_k}{k+1} \sum_{j=1}^k \left[y_i-y^{(i)}_j+B^{(i)}_j\right] + z_{1-\alpha/2} \Delta_\lambda \right], 
\end{equation}
covers $u_i^\star$ with probability at least $100(1-\alpha)\%$.
The width of the plausibility region is given by
$$ 2\sum_{k=1}^{n-1}\frac{\lambda_k}{k+1} \sum_{j=1}^k B^{(i)}_j + 2z_{1-\alpha/2} 
\Delta_\lambda. $$
The optimal set of mixing proportions $\lambda_1,\lambda_2,\cdots,\lambda_{n-1}$ are the solutions to the constrained optimization problem
\begin{align}\label{eq:nPMixOpt}
\min_{\lambda_1,\lambda_2,\cdots,\lambda_{n-1}}\quad & \sum_{k=1}^{n-1}\frac{\lambda_k}{k+1} \sum_{j=1}^k B^{(i)}_j + z_{1-\alpha/2}
\Delta_\lambda \\
\nonumber \text{subject to} \quad & \lambda_1,\lambda_2,\cdots,\lambda_{n-1}\ge0;\ \sum_{k=1}^{n-1}\lambda_k\le1.
\end{align}
We solve problem \eqref{eq:nPMixOpt} numerically and obtain the optimal mixing proportions. Figures 
\ref{fig:n-01},
\ref{fig:n-02},
and 
\ref{fig:n-03}
show numerical results for various configurations of $n$, $M$, and $\gamma$: the blue crosses denote observations from the underlying function $\vartheta_0(t)$ (dashed red line); and the plausibility intervals are displayed alongside.

\ifthenelse{1=0}{}{
\subsection{Asymptotic studies of convergence rate of interval lengths}
\replace{}{
To investigate the asymptotic behavior of our algorithm in terms of both statistical efficiency and computational time, we conducted a sequence of experiments using the simple assertion $A(t) = \{ \vartheta_0: [0,1] \rightarrow \mathbb{R}\ |\ \vartheta_0(t) = 0 \} \subseteq \Theta_{1,1/2}$ with equally spaced observations.

The results of the experiments are summarized in Table \ref{tb:asy} and plotted in Figure \ref{fig:asy}.
We observe that the median widths of the plausibility intervals obtained using our algorithm and the algorithm's runtime exhibit clearly anticipated trends as $n$ increases.
In order to model these two relationships, we apply the logarithm transformation and fit two linear models on the log-scale.

For {\sf median width}, the fitted regression equation is
$$ \log(\mathsf{median\ width}) = 1.378 - 0.267\log(n), $$
with $R^2 = 99.97\%$.
We empirically conclude that the widths of the plausibility intervals decrease at a rate of $\mathcal{O}(n^{-0.267})$.
It is interesting to see that this simulation-based study shows that our empirical asymptotic convergence rate of plausibility intervals for the case of $\gamma=\frac{1}{2}$ is slightly better than or at least close to that in \eqref{eq:asym-rate}, which in this case is $O(n^{-0.250})$.  Although one of these two convergence rates is on point estimation and the other on interval length, such a comparison is arguably meaningful to some extent because IM plausibility intervals are valid in terms of frequency calibration.

For {\sf elapsed time}, the fitted regression equation is
$$ \log(\mathsf{median\ width}) = -10.42 + 4.01\log(n), $$
with $R^2 = 98.67\%$.
We empirically conclude that our algorithm exhibits $\mathcal{O}(n^4)$ time complexity.

\begin{table}[!htb]
  \centering
  \begin{tabular}{rcr}
	\toprule
	$n$ & {\sf Median width} & {\sf Elapsed time} \\
	\midrule
	5	&	2.5941	&	0.04 \\
	10	&	2.1482	&	0.37 \\
	15	&	1.9283	&	1.26 \\
	20	&	1.7772	&	3.19 \\
	25	&	1.6695	&	7.41 \\
	30	&	1.5930	&	16.38 \\
	35	&	1.5307	&	32.40 \\
	40	&	1.4756	&	58.26 \\
	45	&	1.4290	&	105.70 \\
	50	&	1.3910	&	186.41 \\
	60	&	1.3273	&	624.11 \\
	70	&	1.2756	&	1276.15 \\
	80	&	1.2331	&	2014.44 \\
	90	&	1.1976	&	2441.45 \\
	100	&	1.1671	&	3169.96 \\
	\bottomrule
  \end{tabular}
  \caption{Simulation results using the assertion $A(t) = \{ \vartheta_0: [0,1] \rightarrow \mathbb{R}\ |\ \vartheta_0(t) = 0 \} \subseteq \Theta_{1,1/2}$ with equally-spaced observations.}
  \label{tb:asy}
\end{table}

\begin{figure}[!htb]
  \centering
  \subfigure[Plot of {\sf median width} vs $n$]{
    \includegraphics[width=0.48\columnwidth]{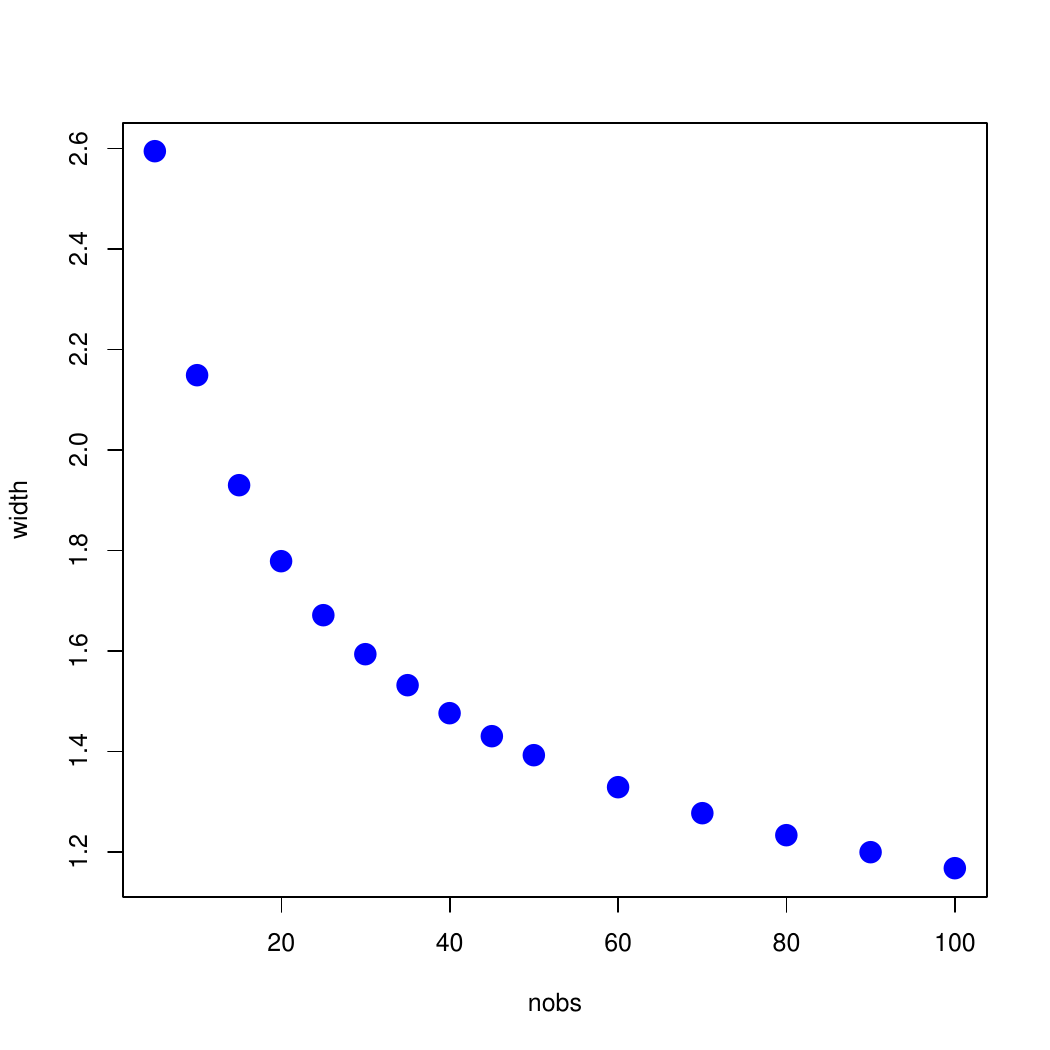}}
  \subfigure[Plot of $\log(\mathsf{median\ width})$ vs $\log(n)$]{
    \includegraphics[width=0.48\columnwidth]{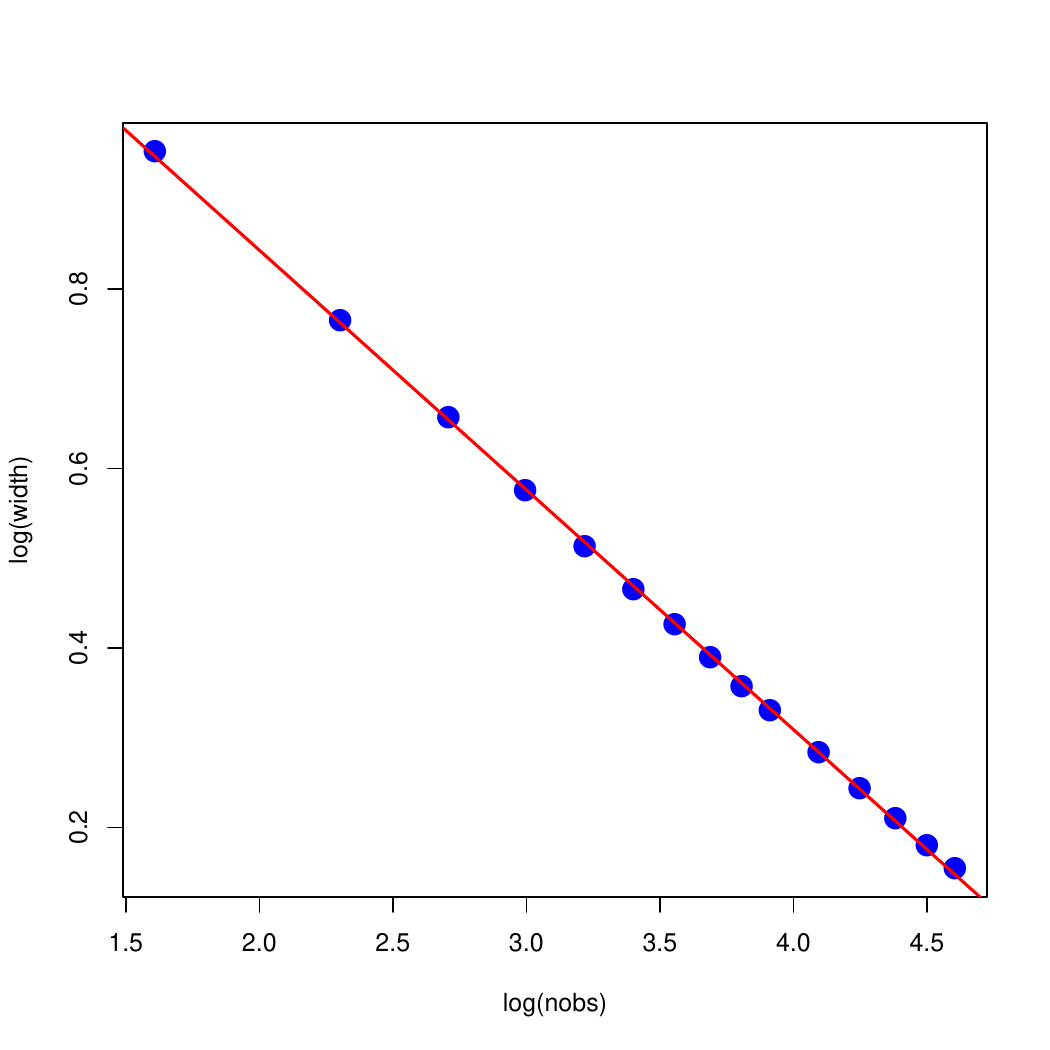}}
  \caption{Visualization of the simulation results in Table \ref{tb:asy}.}
  \label{fig:asy}
\end{figure}
}}

\ifthenelse{1=1}{
\begin{figure}[!htb]
  \centering
  \subfigure[$n=10$]{
    \includegraphics[width=0.48\columnwidth]{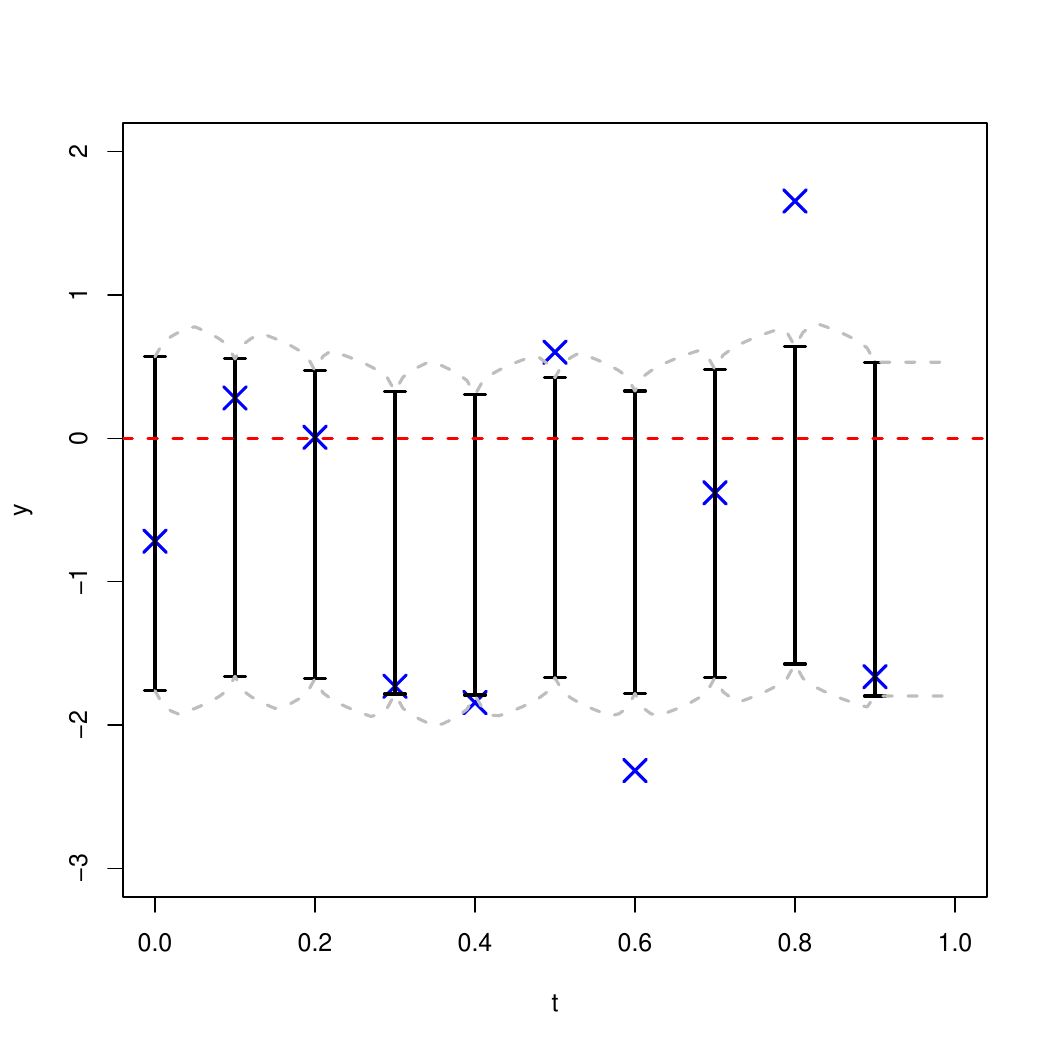}} 
  \subfigure[$n=100$]{
    \includegraphics[width=0.48\columnwidth]{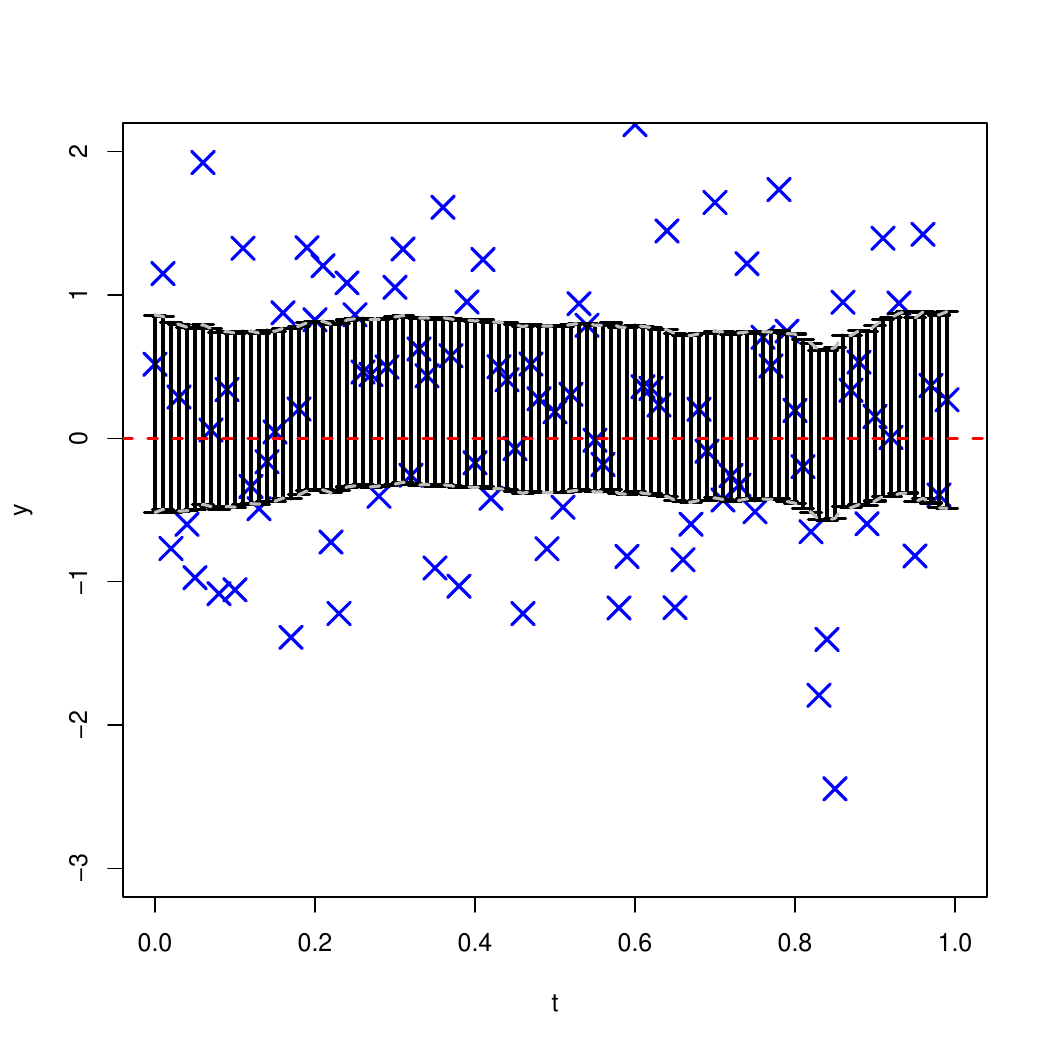}}
  \caption{Visualization of the plausibility intervals constructed using the mixture distribution for the assertion $A(t) = \{ \vartheta_0: [0,1] \rightarrow \mathbb{R}\ |\ \vartheta_0(t) = 0 \} \subseteq \Theta_{1,1/2}$ for various number of observations.}
	\label{fig:n-03}
\end{figure}
}{
\begin{figure}[!htb]
  \centering
  \subfigure[$n=5$]{
    \includegraphics[width=0.48\columnwidth]{asy5}}
  \subfigure[$n=10$]{
    \includegraphics[width=0.48\columnwidth]{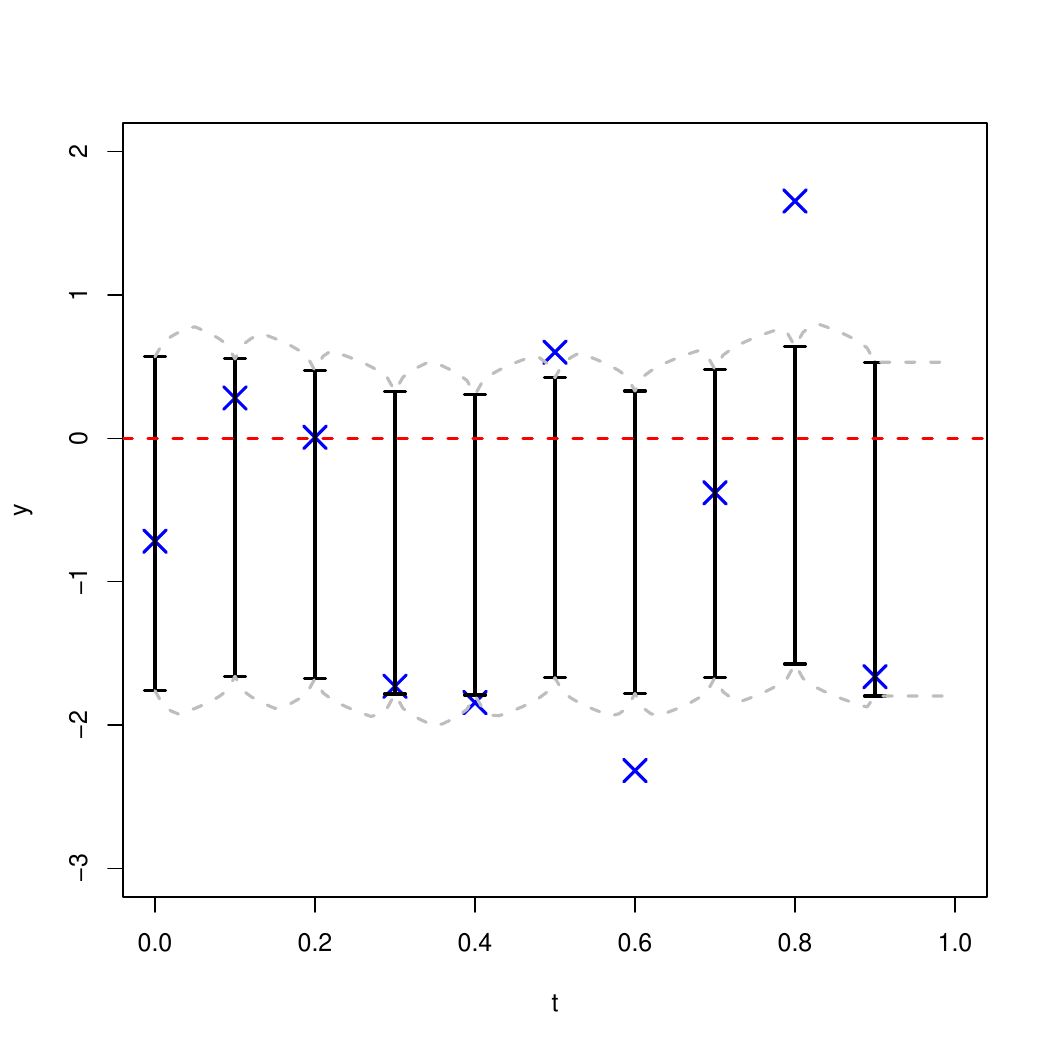}} \\
  \subfigure[$n=15$]{
    \includegraphics[width=0.48\columnwidth]{asy15}}
  \subfigure[$n=20$]{
    \includegraphics[width=0.48\columnwidth]{asy20}}
  \caption{Visualization of the plausibility intervals constructed using the mixture distribution for the assertion $A(t) = \{ \vartheta_0: [0,1] \rightarrow \mathbb{R}\ |\ \vartheta_0(t) = 0 \} \subseteq \Theta_{1,1/2}$ for various number of observations.}
\end{figure}

\begin{figure}[!htb]
  \centering
  \subfigure[$n=25$]{
    \includegraphics[width=0.48\columnwidth]{asy25}}
  \subfigure[$n=30$]{
    \includegraphics[width=0.48\columnwidth]{asy30}} \\
  \subfigure[$n=35$]{
    \includegraphics[width=0.48\columnwidth]{asy35}}
  \subfigure[$n=40$]{
    \includegraphics[width=0.48\columnwidth]{asy40}}
  \caption{Visualization of the plausibility intervals constructed using the mixture distribution for the assertion $A(t) = \{ \vartheta_0: [0,1] \rightarrow \mathbb{R}\ |\ \vartheta_0(t) = 0 \} \subseteq \Theta_{1,1/2}$ for various number of observations.}
\end{figure}

\begin{figure}[!htb]
  \centering
  \subfigure[$n=45$]{
    \includegraphics[width=0.48\columnwidth]{asy45}}
  \subfigure[$n=50$]{
    \includegraphics[width=0.48\columnwidth]{asy50}} \\
  \subfigure[$n=60$]{
    \includegraphics[width=0.48\columnwidth]{asy60}}
  \subfigure[$n=70$]{
    \includegraphics[width=0.48\columnwidth]{asy70}}
  \caption{Visualization of the plausibility intervals constructed using the mixture distribution for the assertion $A(t) = \{ \vartheta_0: [0,1] \rightarrow \mathbb{R}\ |\ \vartheta_0(t) = 0 \} \subseteq \Theta_{1,1/2}$ for various number of observations.}
\end{figure}

\begin{figure}[!htb]
  \centering
  \subfigure[$n=80$]{
    \includegraphics[width=0.48\columnwidth]{asy80}}
  \subfigure[$n=90$]{
    \includegraphics[width=0.48\columnwidth]{asy90}} \\
  \subfigure[$n=100$]{
    \includegraphics[width=0.48\columnwidth]{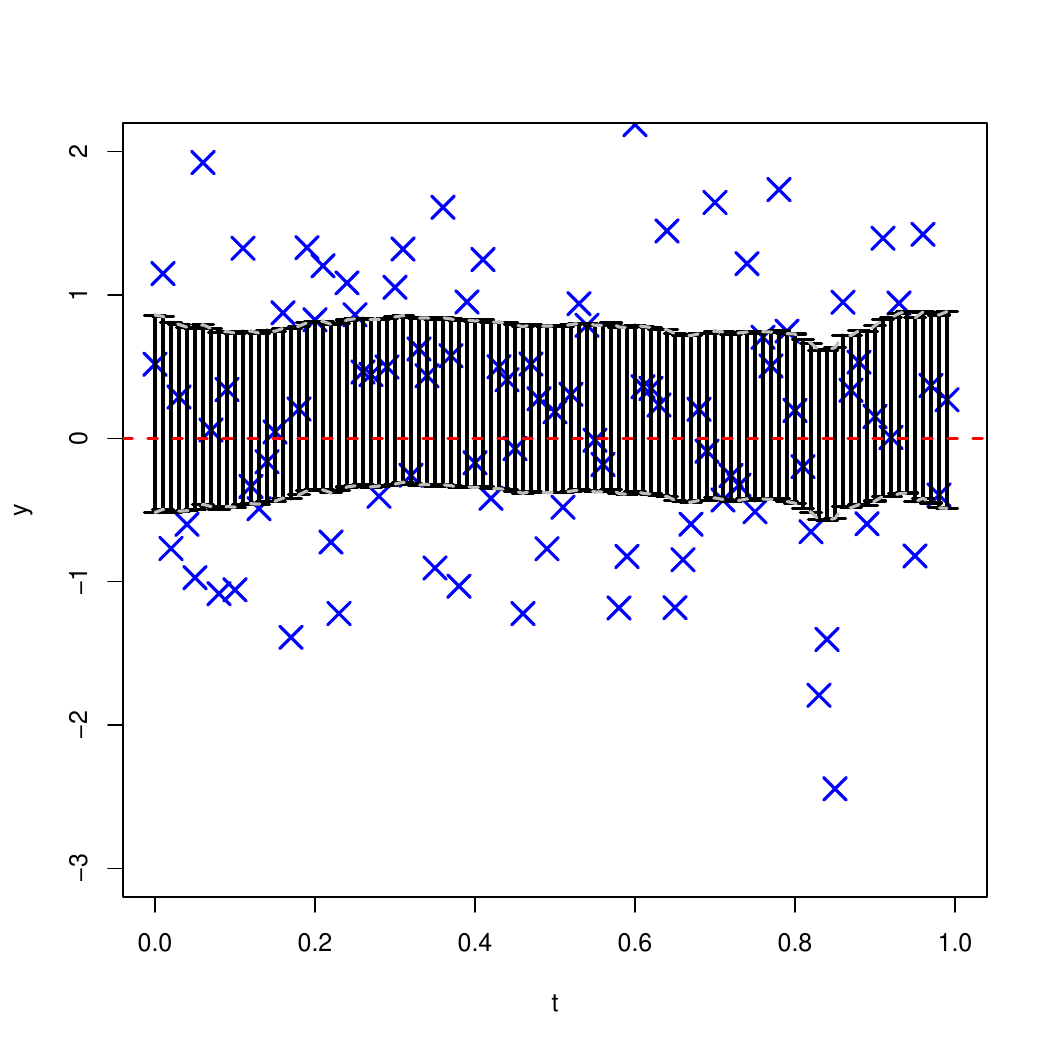}}
  \caption{Visualization of the plausibility intervals constructed using the mixture distribution for the assertion $A(t) = \{ \vartheta_0: [0,1] \rightarrow \mathbb{R}\ |\ \vartheta_0(t) = 0 \} \subseteq \Theta_{1,1/2}$ for various number of observations.}
\end{figure}
}


\ifthenelse{1=1}{}{
\begin{lstlisting}
########## The n-Point Case ##############

set.seed(1234)

order_dist <- function(t, i, arg = TRUE, discard = TRUE) {
  # Sorts the elements of a vector in ascending distance to a specific element.
  #
  # Args:
  #   t: A vector whose elements are to be sorted.
  #   i: An index specifying the element to which distances are computed.
  #   arg: If TRUE, returns the sorted indices.
  #        If FALSE, returns the sorted elements. Default is TRUE.
  #   discard: If TRUE, discards the first index or element in the sorted vector,
  #               which would always be i or t[i] itself.
  # If FALSE, keep the first index or element in the sorted vector.
  #
  # Returns:
  #   The sorted indices or elements of t in ascending distance to t[i].
  if (discard == TRUE) {
    if (arg == TRUE) {
      order(abs(t - t[i]))[2:length(t)]
    } else {
      sort(abs(t - t[i]))[2:length(t)]
    }
  }
  else {
    if (arg == TRUE)
      order(abs(t - t[i]))
    else
      sort(abs(t - t[i]))
  } 
}

npregpl <- function(data, M = 1, gam = 1/2, sigma = 1, alpha = 0.05) {
  # Computes point-wise plausibility limits for nonparametric regression.
  #
  # Args:
  #   data: An n*2 matrix specifying the n pairs of observations, whose
  #         first column records the design points, and whose
  #         second column records the observed response (NA denotes missing).
  #   M, gam: Holder space parameters. The default is M = 1, gam = 1/2.
  #   sigma: Standard deviation. Default is s = 1.
  #   alpha: Significance level. Default is alpha = 0.05.
  #
  # Returns:
  #   X: An n*4 matrix whose first two columns consist of the original data, and
  #      whose last two columns specify the lower and upper plausibility limits.
  #      The rows of X are sorted in ascending order of the design points.
  N <- dim(data)[1]  # Total number of points
  X <- cbind(data[order(data[, 1]), ], matrix(rep(0, N * 2), N, 2))
  obs <- (1:N)[is.na(X[, 2]) == 0]  # Observed indices
  mis <- (1:N)[is.na(X[, 2]) == 1]  # Missing indices
  n <- length(obs)
  if (n == 0) {
    print("Error: no observed points!")
    return(NA)
  }
  t <- X[obs, 1]  # Design points
  y <- X[obs, 2] / sigma  # Normalized response
  colnames(X) <- c("t", "y", "lower", "upper")
  z <- qnorm(1 - alpha / 2)  # Normal quantile
  for (i in 1:n) {
    # Predict the i-th point
    t_i <- t[order_dist(t, i)]
    y_i <- y[order_dist(t, i)]
    B_i <- M * abs(t_i - t[i]) ^ gam
    mixture <- function(B_i) {
      # Determine the optimal mixture distribution
      mixture_var <- function(lambda) {
        var <- (1 - sum(lambda * (1:(n - 1)) / (2:n))) ^ 2
        for (j in 1:(n - 1))
          var = var + sum(lambda[j:(n - 1)] / ((j + 1):n)) ^ 2
        var
      }
      mixture_width <- function(lambda) {
        pivot <- 0
        for (k in 1:(n - 1))
          pivot = pivot + 2 * lambda[k] * sum(B_i[1:k]) / (k+1)
        pivot + 2 * z * sqrt(mixture_var(lambda))
      }
      mixture_width_grad <- function(lambda) {
        grad <- rep(0, n - 1)
        for (k in 1:(n - 1)) {
          var_grad_k <- 2 * (1 - sum(lambda * (1:(n - 1)) / (2:n))) *
                            (-k / (k + 1))
          for (j in 1:k)
              var_grad_k = var_grad_k +
                           2 * sum(lambda[j:(n - 1)]/((j + 1):n)) / (k + 1)
          grad[k] = 2 * sum(B_i[1:k]) / (k + 1) +
                    z / sqrt(mixture_var(lambda)) * var_grad_k
        }
        grad
      }
      lambda_init <- runif(n - 1, 0, 1 / (n - 1))  # Random initialization
      constrOptim(theta = lambda_init,
                  f = mixture_width,
                  grad = mixture_width_grad,
                  ui = rbind(diag(n - 1), rep(-1, n - 1)),
                  ci = c(rep(0, n - 1), -1))
    }
    mixture_opt <- mixture(B_i)
    lambda <- mixture_opt$par
    width <- mixture_opt$value
    center <- 0
    for (k in 1:(n - 1))
      center = center + lambda[k] * sum(y[i] - y_i[1:k]) / (k + 1)
    upper <- center + width / 2
    lower <- center - width / 2
    X[obs[i], 3:4] <- sigma * (y[i] - c(upper, lower))
  }
  for (i in 1:length(mis)) {
    ti <- X[mis[i], 1]
    l <- sum(obs - mis[i] < 0)
    r <- sum(obs - mis[i] > 0)
    if (l == 0) {
      tr <- X[obs[1], 1]
      Br <- M * abs(ti - tr) ^ gam
      X[mis[i], 3] = X[obs[1], 3] - Br
      X[mis[i], 4] = X[obs[1], 4] + Br
    } else if (r == 0) {
      tl <- X[obs[n], 1]
      Bl <- M * abs(ti - tl) ^ gam
      X[mis[i], 3] = X[obs[n], 3] - Br
      X[mis[i], 4] = X[obs[n], 4] + Br
    } else {
      tl <- X[obs[l], 1]
      tr <- X[obs[n - r + 1], 1]
      Bl <- M * abs(ti - tl) ^ gam
      Br <- M * abs(ti - tr) ^ gam
      X[mis[i], 3] = max(X[obs[l], 3] - Bl, X[obs[n - r + 1], 3] - Br)
      X[mis[i], 4] = min(X[obs[l], 4] + Bl, X[obs[n - r + 1], 4] + Br)
    }
  }
  X
}

##### Example: the square root function #####

require(plotrix)

sqrt_pl <- function(nobs, sigma) {
  nmis <- 100
  tmis <- seq(0, 1 - 1 / nmis, 1 / nmis)
  ymis <- rep(NA, nmis)
  tobs <- runif(nobs)
  yobs <- sqrt(tobs) + rnorm(nobs, sd = sigma)
  data <- cbind(c(tmis, tobs), c(ymis, yobs))
  X <- npregpl(data, sigma = sigma)
  N <- dim(data)[1]
  obs <- (1:N)[is.na(X[, 2]) == 0]  # Observed indices
  mis <- (1:N)[is.na(X[, 2]) == 1]  # Missing indices
  plot(X[obs, 1], X[obs, 2],
       xlim = c(0, 1), ylim = c(-2, 3),
       xlab = "t", ylab = "y",
       pch = 4, cex = 2, lwd = 2, col = "blue")
  lines(X[, 1], X[, 3], lty = "dashed", lwd = 2, col = "gray")
  lines(X[, 1], X[, 4], lty = "dashed", lwd = 2, col = "gray")
  plotCI(x = X[obs, 1],
         y = (X[obs, 4] + X[obs, 3]) / 2,
         uiw = (X[obs, 4] - X[obs, 3]) / 2,
         xlab = "t", ylab = "y",
         pch = NA, lwd = 2, add = TRUE)
  curve(sqrt, 0, 1, lty = "dashed", lwd = 2, col = "red", add = TRUE)
}

# Experiments
NOBS <- c(5, 10, 20)
SIGMA <- c(1, 1/2)
for (i in 1:length(NOBS)) {
  nobs <- NOBS[i]
  for (j in 1:length(SIGMA)) {
    sigma <- SIGMA[j]
    pdf(paste("sqrt", i, j, ".pdf", sep = ""))
    sqrt_pl(nobs, sigma)
    dev.off()    
  }
}

##### Example: the zero function #####

zero_pl <- function(nobs, M, gam, sigma = 1) {
  nmis <- 100
  tmis <- seq(0, 1 - 1 / nmis, 1 / nmis)
  ymis <- rep(NA, nmis)
  tobs <- runif(nobs)
  yobs <- rnorm(nobs, sd = sigma)
  data <- cbind(c(tmis, tobs), c(ymis, yobs))
  X <- npregpl(data, M = M, gam = gam, sigma = sigma)
  N <- dim(data)[1]
  obs <- (1:N)[is.na(X[, 2]) == 0]  # Observed indices
  mis <- (1:N)[is.na(X[, 2]) == 1]  # Missing indices
  plot(X[obs, 1], X[obs, 2],
       xlim = c(0, 1), ylim = c(-3, 2),
       xlab = "t", ylab = "y",
       pch = 4, cex = 2, lwd = 2, col = "blue")
  lines(X[, 1], X[, 3], lty = "dashed", lwd = 2, col = "gray")
  lines(X[, 1], X[, 4], lty = "dashed", lwd = 2, col = "gray")
  plotCI(x = X[obs, 1],
         y = (X[obs, 4] + X[obs, 3]) / 2,
         uiw = (X[obs, 4] - X[obs, 3]) / 2,
         xlab = "t", ylab = "y",
         pch = NA, lwd = 2, add = TRUE)
  abline(a = 0, b = 0, lty = "dashed", lwd = 2, col = "red")
}

# Experiments
Ms <- c(1, 1/2)
GAMs <- c(1/5, 1)
for (i in 1:length(Ms)) {
  M <- Ms[i]
  for (j in 1:length(GAMs)) {
    gam <- GAMs[j]
    pdf(paste("zero", i, j, ".pdf", sep = ""))
    zero_pl(nobs = 15, M = M, gam = gam)
    dev.off()
  }
}

##### Asymptotic studies using the zero function #####

asy_pl <- function(nobs, sigma = 1) {
  nmis <- 100
  tmis <- seq(0, 1 - 1 / nmis, 1 / nmis)
  ymis <- rep(NA, nmis)
  tobs <- seq(0, 1 - 1 / nobs, 1 / nobs)  # Equally-spaced observations
  yobs <- rnorm(nobs, sd = sigma)
  data <- cbind(c(tmis, tobs), c(ymis, yobs))
  X <- npregpl(data, sigma = sigma)
  N <- dim(data)[1]
  obs <- (1:N)[is.na(X[, 2]) == 0]  # Observed indices
  mis <- (1:N)[is.na(X[, 2]) == 1]  # Missing indices
  plot(X[obs, 1], X[obs, 2],
       xlim = c(0, 1), ylim = c(-3, 2),
       xlab = "t", ylab = "y",
       pch = 4, cex = 2, lwd = 2, col = "blue")
  lines(X[, 1], X[, 3], lty = "dashed", lwd = 2, col = "gray")
  lines(X[, 1], X[, 4], lty = "dashed", lwd = 2, col = "gray")
  plotCI(x = X[obs, 1],
         y = (X[obs, 4] + X[obs, 3]) / 2,
         uiw = (X[obs, 4] - X[obs, 3]) / 2,
         xlab = "t", ylab = "y",
         pch = NA, lwd = 2, add = TRUE)
  abline(a = 0, b = 0, lty = "dashed", lwd = 2, col = "red")
  median(X[obs, 4] - X[obs, 3])  # Median width of all plausibility intervals
}

# Experiments

trials <- c(seq(5, 50, 5), seq(60, 100, 10))
ntrials <- length(trials)
results <- matrix(rep(0, ntrials * 3), ntrials, 3)
colnames(results) <- c("nobs", "width", "time")
for (i in 1:ntrials) {
  nobs <- trials[i]
  pdf(paste("asy", nobs, ".pdf", sep = ""))
  begin <- proc.time()
  width <- asy_pl(nobs)
  end <- proc.time()
  dev.off()
  results[i, ] <- c(nobs, width, (end - begin)[3])
}
write.table(results, "results.txt", sep = "\t")

nobs <- results[, 1]
width <- results[, 2]
time <- results[, 3]

pdf("asy.pdf")
plot(nobs, width, pch = 16, cex = 2, col = "blue")
dev.off()

fit_asy <- lm(log(width) ~ log(nobs))  # Fit linear model on log-scale
summary(fit_asy)
pdf("asylog.pdf")
plot(log(nobs), log(width), pch = 16, cex = 2, col = "blue")
abline(fit_asy, lwd = 2, col = "red")
dev.off()

pdf("time.pdf")
plot(nobs, time, pch = 16, cex = 2, col = "blue")
dev.off()

fit_time <- lm(log(time) ~ log(nobs))  # Fit linear model on log-scale
summary(fit_time)
pdf("timelog.pdf")
plot(log(nobs), log(time), pch = 16, cex = 2, col = "blue")
abline(fit_time, lwd = 2, col = "red")
dev.off()
\end{lstlisting}
}

\section{Discussion}
\label{s:Discussion}

\replace{}{
We have developed a \emph{partial conditioning} method to extend IMs 
for producing valid and efficient inference about
many-normal-means when the means are subject to H\"older constraints. The problem was motivated by the challenging setting in nonparametric regression, where no known methods are available to produce valid confidence intervals, even asymptotically. Thus, we make no attempt to compare the proposed method to existing methods, as our focus here is to extend conditional IMs to tackle many-normal-means with H\"older constraints (including the Lipschitz condition as a special case).

Nevertheless, 
the simple simulation-based study in Section \ref{s:GeneralCase} shows that our empirical asymptotic convergence rate of plausibility intervals for the case of $\gamma=\frac{1}{2}$ is about $O(n^{-0.267})$. This is slightly better than or at least close to that in \eqref{eq:asym-rate}, which in this case is $O(n^{-0.250})$.  As augued in Section \ref{s:GeneralCase}, such a comparison is arguably meaningful. Perhaps, the implication is of extreme importance to tackling the challenging problem of constructing confidence interval in nonparametric regression, which has been perceived as unsovable with existing approaches.

In addition, we expect the proposed method to inspire novel applications and theoretical developments, such as in nonparametric regression \citep[c.f.][and references therein]{lepskii1991problem,wang2013uniform}.
In the case when point estimation is of interest, which is often the starting-point for frequentist methods, our results \eqref{eq:plint-n} for the general $n$ case can be used to provide an alternative local shrinkage or smoothing scheme, in the same manner as described in Section~\ref{ss:skrinkage}.

Furthermore, while we have focused on developing the method for 
many-normal-means in this work, it can be extended to handle
many-binomial-means and many-Poisson-means with H\"older constraints.
Nevertheless, a more general formulation of partial conditioning and its mathematical theory deserve future development.
}

\section*{Acknowledgements}
\replace{}{The authors are grateful to the Editor and two referees for their insightful comments and constructive suggestions that have led to significant improvement of this article.}
\ifthenelse{1=1}{}{

In this section, we discuss methods for estimating the variance parameter $\sigma^2$.


\subsubsection{Conditional Distributions}

Let us first consider the case of $n=3$, where we have three pairs of observations $(t_i, Y_i), i=1,2,3$. Without loss of generality, we assume that $t_1\le t_2\le t_3$. For $i = 2,3$, we define the differences
$$ \Delta_i := Y_i - Y_{i-1} = [\vartheta_0(t_i) - \vartheta_0(t_{i-1})] + [U_i - U_{i-1}],\ i = 2, 3. $$
The joint distribution of $\Delta_2$ and $\Delta_3$ is given by
\begin{equation}
\left[ \begin{array}{c} \Delta_2 \\ \Delta_3 \end{array} \right] \sim\mathcal{N}_2 \left( \bm{\mu} = \left[\begin{array}{c} d_2 \\ d_3 \end{array}\right],\ \bm{\Sigma} = \sigma^2 \left[\begin{array}{ccc} 2 & -1 \\ -1 & 2 \end{array}\right] \right),
\end{equation}
where we define $d_i := \vartheta_0(t_i) - \vartheta_0(t_{i-1}),\ i= 2,3$.
Multiplying both sides of the equation by
$$ \bm{\Sigma}^{-1/2} = \frac{1}{\sqrt{12 + 6\sqrt{3}}} \left[ \begin{array}{cc} 2 + \sqrt{3} & 1 \\ 1 & 2 + \sqrt{3} \end{array} \right]  $$
yields
\begin{equation}
\bm{\Sigma}^{-1/2} \left[ \begin{array}{c} \Delta_2 - d_2 \\ \Delta_3 - d_3 \end{array} \right] \sim\mathcal{N}_2 \left( \bm{0}, \sigma^2 \bm{I}\right),
\end{equation}
which is equivalent to
\begin{equation}\label{eq:Var3P}
\begin{dcases}
\tilde{\Delta}_2 = \tilde{d}_2 + \sigma Z_2; \\
\tilde{\Delta}_3 = \tilde{d}_3 + \sigma Z_3, \\
\end{dcases}
\end{equation}
where
\begin{align*}
\tilde{\Delta}_2 & = [(2 + \sqrt{3})\,\Delta_2 + \Delta_3] / \sqrt{12 + 6\sqrt{3}}; \\
\tilde{\Delta}_3 & = [\Delta_2 + (2 + \sqrt{3})\,\Delta_3] / \sqrt{12 + 6\sqrt{3}}, \\
\tilde{d}_2 & = [(2 + \sqrt{3})\,d_2 + d_3] / \sqrt{12 + 6\sqrt{3}}; \\
\tilde{d}_3 & = [d_2 + (2 + \sqrt{3})\,d_3] / \sqrt{12 + 6\sqrt{3}},
\end{align*}
and $Z_2, Z_3\ i.i.d. \sim \mathcal{N}(0,1)$.

Taking absolute values and then logarithms on both sides of the equations in \eqref{eq:Var3P}, we obtain:
\begin{equation}\label{eq:logVar3P}
\begin{dcases}
\log|\tilde{\Delta}_2 - \tilde{d}_2| = \log\sigma + \log|Z_2|; \\
\log|\tilde{\Delta}_3 - \tilde{d}_3| = \log\sigma + \log|Z_3|. \\
\end{dcases}
\end{equation}

Notice that equation \eqref{eq:logVar3P} implies that
$$ \log|Z_3| - \log|Z_2| = \log|\tilde{\Delta}_3 - \tilde{d}_3| - \log|\tilde{\Delta}_2 - \tilde{d}_2|. $$

For convenience, we introduce the notations $U_i = \log|Z_i|, i = 2,3$ and $V = U_3 - U_2 = \log|Z_3/Z_2|$.
To predict $U_3$, we could utilize either the conditional distribution
\begin{equation*}
U_3\,\bigg| \left\{V = \log|\tilde{\Delta}_3 - \tilde{d}_3| - \log|\tilde{\Delta}_2 - \tilde{d}_2|\right\}\ \sim\ f^C_{U_3|V}(u|v),
\end{equation*}
or the marginal distribution
\begin{equation}\label{eq:Var3PMarg}
U_3 \sim f^M(u) = \sqrt{\frac{2}{\pi}} \exp\left\{u - \frac{1}{2} e^{2u}\right\},
\end{equation}
which follows from the fact that $Z_3 \sim \mathcal{N}(0,1)$.

Let us now derive the form of the conditional distribution.
We have $U_2, U_3\ i.i.d. \sim f^M$, and the transformation:
\begin{equation*}
\begin{dcases}
U_3 = U_3; \\
V = U_3 - U_2.
\end{dcases}
\end{equation*}
The joint density of $(U_3,V)$ is given by
\begin{align*}
f_{U_3,\,V}(u, v)
& = f_{U_3,\,U_2}(u, u-v) = f_{U_3}(u)\,f_{U_2}(u-v) \\
& = \frac{2}{\pi} \exp\left\{(2u-v) - \frac{1}{2}e^{2u}\left(1 + e^{-2v}\right) \right\}, \quad u,v\in(-\infty,\infty);
\end{align*}
and the marginal density of $V$ is given by
\begin{align*}
f_V(v)
& = \int_{-\infty}^\infty f_{U_3,\,V}(u, v) \D u \\
& = \frac{2}{\pi} \int_{-\infty}^\infty \exp\left\{(2u-v) - \frac{1}{2}e^{2u}\left(1 + e^{-2v}\right) \right\} \D u \\
& = \frac{1}{\pi} e^{-v} \int_{-\infty}^\infty \exp\left\{2u - \frac{1}{2}e^{2u}\left(1 + e^{-2v}\right) \right\} \D u \\
& = \frac{1}{\pi} e^{-v} \left[-\frac{2\exp\{-\frac{1}{2}(1 + e^{-2v})\,e^{2u}\}}{(1 + e^{-2v})}\right] \bigg|_{-\infty}^\infty \\
& = \frac{2e^{-v}}{\pi(1 + e^{-2v})} = \frac{\text{sech}(v)}{\pi},
\end{align*}
which we recognize as the \emph{hyperbolic secant distribution}.
Therefore, the conditional density
\begin{align}
f^C_{U_3|V}(u)
& \nonumber = \frac{f_{U_3,\,V}(u, v)}{f_V(v)} \\
& = \left(1 + e^{-2v}\right)\exp\left\{2u - \frac{1}{2} \left(1 + e^{-2v}\right) e^{2u} \right\}, \quad u\in(-\infty,\infty).
\end{align}
Notice that $e^{U_3} | V = v$ follows a $\mathsf{Weibull}\, (\mathsf{scale}=\sqrt{2/(1+e^{-2v})},\,\mathsf{shape}=2)$ density symmetrically extended on $(-\infty, \infty)$.

\subsubsection{Prediction via Mixture Distributions}

In order to predict $U_3$, we utilize a mixture of the conditional density and the marginal density, which results in the mixture density:
\begin{equation}
f_{U_3}(u; \lambda, v) = \lambda \left(1 + e^{-2v}\right)\exp\left\{2u - \frac{1}{2} \left(1 + e^{-2v}\right) e^{2u} \right\} + (1-\lambda)\sqrt{\frac{2}{\pi}} \exp\left\{u - \frac{1}{2} e^{2u}\right\},\ u\in(-\infty,\infty),
\end{equation}
and whose CDF is given by
\begin{equation}\label{eq:Var3PMixCDF}
F_{U_3}(u;\lambda,v) = \lambda \left[ 1 - \exp\left\{-\frac{1}{2} \left(1 + e^{-2v}\right) e^{2u} \right\} \right] + (1-\lambda) \left[ 2\Phi(e^u) - 1 \right], \quad u\in(-\infty,\infty),
\end{equation}
where
\begin{align}
v
\nonumber & = \log|\tilde{\Delta}_3 - \tilde{d}_3| -  \log|\tilde{\Delta}_2 - \tilde{d}_2| \\
\nonumber & = \log \left|(\Delta_2-d_2) + (2+\sqrt{3})(\Delta_3-d_3) \right| - \log \left|(2+\sqrt{3})(\Delta_2-d_2) + (\Delta_3-d_3) \right|  \\
\nonumber & = \log\left| \left[\Delta_2 + (2+\sqrt{3})\Delta_3 \right] - \left[d_2 + (2+\sqrt{3}) d_3 \right]\right| - \log\left| \left[(2+\sqrt{3})\Delta_2 + \Delta_3\right] - \left[(2+\sqrt{3}) d_2 + d_3 \right]\right|. \\
& \in \left[ \log \left| \frac{\left|(2+\sqrt{3})\Delta_2 + \Delta_3\right| - \left[(2+\sqrt{3}) B_{21} + B_{32}\right]}{\Delta_2 + (2+\sqrt{3}) \Delta_3 + B_{21} + (2+\sqrt{3}) B_{32}} \right|, \ \log \left| \frac{\Delta_2 + (2+\sqrt{3}) \Delta_3 + B_{21} + (2+\sqrt{3}) B_{32}}{\left|(2+\sqrt{3})\Delta_2 + \Delta_3\right| - \left[(2+\sqrt{3}) B_{21} + B_{32}\right]} \right| \, \right], \label{eq:Var3PvInt}
\end{align}
where we recall that $\Delta_i = Y_i - Y_{i-1}$ are known constants, and $d_i = \vartheta(t_i) - \vartheta(t_{i-1}) \in [-B_{i,\,i-1}, B_{i,\,i-1}]$.

For the marginal case ($\lambda = 0$), setting $F_{U_3}(u; 0, v)$ to be $\alpha/2$ and $1 - \alpha/2$ respectively and then solving for $u$ yields the $(1-\alpha)\%$ plausibility interval\footnote{Whether this method of constructing plausibility intervals is optimal remains to be examined.}:
$$ \left[ \log\left\{ \Phi^{-1} \left(\frac{\alpha}{4} + \frac{1}{2} \right) \right\},\ \log\left\{ \Phi^{-1} \left( 1 - \frac{\alpha}{4} \right) \right\} \right]. $$
For $\alpha=0.05$, the length of the 95\% plausibility interval is equal to 4.2700.

For the conditional case ($\lambda = 1$), setting $F_{U_3}(u; 1, v)$ to be $\alpha/2$ and $1 - \alpha/2$ respectively and then solving for $u$ yields the $(1-\alpha)\%$ plausibility interval for each fixed value of $v$:
$$ \left[ \frac{1}{2} \log \left\{\frac{-2\log(1-\alpha/2)}{1+e^{-2v}} \right\},\ \frac{1}{2} \log \left\{\frac{-2\log(\alpha/2)}{1+e^{-2v}} \right\} \right], $$
whose length can be calculated as 
$$ \frac{1}{2} \log\left\{ \frac{\log(\alpha/2)}{\log(1-\alpha/2)} \right\}, $$
which, interestingly, does not depend on $v$. For $\alpha=0.05$, the length of the 95\% plausibility interval is equal to 2.4908.
Also note that the plausibility limits are monotonically decreasing in $v$.

Simulation shows that the optimal mixing proportion of the mixture plausibility interval \eqref{eq:Var3PMixCDF} appears to be always equal to 1, regardless of the values of $v$. In this case, the $(1-\alpha)\%$ plausibility interval for $U_3$ is given by
\begin{align}
\bigcup_{v\in\eqref{eq:Var3PvInt}} \left[ \frac{1}{2} \log \left\{\frac{-2\log(1-\alpha/2)}{1+e^{-2v}} \right\},\ \frac{1}{2} \log \left\{\frac{-2\log(\alpha/2)}{1+e^{-2v}} \right\} \right].
\end{align}

Returning to equation \eqref{eq:logVar3P}, we arrive at the 95\% plausibility interval for $\log\sigma$:
\begin{align*}
\log\sigma
& = \log|\tilde{\Delta}_3 - \tilde{d}_3| - \log|Z_3| \\
& \in \left[ \log \left| |(2+\sqrt{3})\Delta_2 + \Delta_3| - \left[(2+\sqrt{3}) B_{21} + B_{32}\right] \right|, \ \log \left| \Delta_2 + (2+\sqrt{3}) \Delta_3 + B_{21} + (2+\sqrt{3}) B_{32} \right| \, \right] \\
& - \left[ -\frac{1}{2} \log \left\{\frac{-2\log(\alpha/2)}{1+e^{-2v_1}} \right\},\ -\frac{1}{2} \log \left\{\frac{-2\log(1-\alpha/2)}{1+e^{-2v_2}} \right\} \right],
\end{align*}
where
\begin{align*}
\begin{dcases}
v_1 = \log \left| \frac{\left|(2+\sqrt{3})\Delta_2 + \Delta_3\right| - \left[(2+\sqrt{3}) B_{21} + B_{32}\right]}{\Delta_2 + (2+\sqrt{3}) \Delta_3 + B_{21} + (2+\sqrt{3}) B_{32}} \right|; \\
v_2 = \log \left| \frac{\Delta_2 + (2+\sqrt{3}) \Delta_3 + B_{21} + (2+\sqrt{3}) B_{32}}{\left|(2+\sqrt{3})\Delta_2 + \Delta_3\right| - \left[(2+\sqrt{3}) B_{21} + B_{32}\right]} \right|.
\end{dcases}
\end{align*}

\newpage
}


\ifthenelse{1=1}{}{
In this section we explore the idea of \emph{partial conditioning}.
We begin by illustrating this rather general concept under the setting of Section \ref{sec:TwoPoint}, where we considered pointwise inference for the two-point case of the nonparametric regression problem.

Recall the discussion in Section \ref{sec:TwoPoint}, where we proposed three different approaches to conduct inference on the auxiliary variable $U_2$ by utilizing the marginal, conditional, or mixture density. Let us now revisit these approaches from a different perspective.
More specifically, instead of working directly at the density-level, we convert the problem into constructing an appropriate probability measure on (nested) predictive random sets which allows us to conduct \emph{valid} and \emph{efficient} inference.


\begin{proposition}\label{prop:TwoPointPC}
The three different approaches proposed in Section \ref{sec:TwoPoint} lead to both valid and efficient inference, that is, using the default nested random sets $\mathcal{S} = [-|\delta|,\,|\delta|]$, we have:
\begin{itemize}
\item Prediction using the marginal distribution $U\sim\mathcal{N}(0, 1)$ is both valid and efficient, i.e., 
$$ \mathsf{P}_\mathcal{S}\{U\in\mathcal{S}\}\sim\mathcal{U}(0, 1). $$
\item  Prediction using the conditional distribution 
$$ U|V = v^\star\sim\mathcal{N}\left(\frac{v^\star}{2},\,\frac{1}{2}\right) $$
is both valid and efficient, i.e.,
$$ \mathsf{P}_\mathcal{S}\left\{U\in\mathcal{S} + \frac{v^\star}{2}\right\}\sim\mathcal{U}(0, 1). $$
\item For any $\lambda\in[0,1]$ fixed, the mixture distribution
$$ U|V = v^\star\sim\mathcal{N}\left(\frac{\lambda}{2}v^\star,\,1 - \lambda + \frac{\lambda^2}{2}\right) $$
is both valid and efficient, i.e., 
$$ \mathsf{P}_\mathcal{S}\left\{U\in\mathcal{S} + \frac{\lambda\,v^\star}{2}\right\}\sim\mathcal{U}(0, 1). $$
\end{itemize}
\end{proposition}

\begin{proof}
Marginal case: for any $\alpha\in(0,1)$, we have that
\begin{align*}
\mathsf{P}_U\left\{\mathsf{P}_\mathcal{S}\left\{U\in[-|\delta|,\,|\delta|]\right\} \le \alpha\right\}
& = \mathsf{P}_U\left\{\mathsf{P}_\mathcal{S}\left\{|\delta| \ge \left|U\right|\right\} \le \alpha\right\} \\
& = \mathsf{P}_U\left\{2\left[1 - \Phi\left(\left|U \right|\right)\right] \le \alpha\right\} \\
& = \mathsf{P}_U\left\{ |U| \ge \Phi^{-1}\left(1-\frac{\alpha}{2}\right) \right\} \\
& = 2\left[1 - \left(1-\frac{\alpha}{2}\right)\right] \\
& = \alpha.
\end{align*}
Conditional case: for any $\alpha\in(0,1)$, we have that
\begin{align*}
\mathsf{P}_U\left\{\mathsf{P}_\mathcal{S}\left\{U\in[-|\delta|,\,|\delta|] + \frac{\,v^\star}{2}\right\} \le \alpha\right\}
& = \mathsf{P}_U\left\{\mathsf{P}_\mathcal{S}\left\{|\delta| \ge \left|U - \frac{v^\star}{2}\right|\right\} \le \alpha\right\} \\
& = \mathsf{P}_U\left\{2\left[1 - \Phi\left(\sqrt{2}\left|U - \frac{v^\star}{2}\right|\right)\right] \le \alpha\right\} \\
& = \mathsf{P}_U\left\{ \sqrt{2}\left|U - \frac{\,v^\star}{2}\right| \ge \Phi^{-1}\left(1-\frac{\alpha}{2}\right) \right\} \\
& = 2\left[1 - \left(1-\frac{\alpha}{2}\right)\right] \\
& = \alpha.
\end{align*}
Mixture case: for any $\alpha\in(0,1)$, we have that
\begin{align*}
\mathsf{P}_U\left\{\mathsf{P}_\mathcal{S}\left\{U\in[-|\delta|,\,|\delta|] + \frac{\lambda\,v^\star}{2}\right\} \le \alpha\right\}
& = \mathsf{P}_U\left\{\mathsf{P}_\mathcal{S}\left\{|\delta| \ge \left|U - \frac{\lambda\,v^\star}{2}\right|\right\} \le \alpha\right\} \\
& = \mathsf{P}_U\left\{2\left[1 - \Phi\left(\frac{\left|U - \frac{\lambda\,v^\star}{2}\right|}{\sqrt{1 - \lambda + \frac{\lambda^2}{2}}}\right)\right] \le \alpha\right\} \\
& = \mathsf{P}_U\left\{ \frac{\left|U - \frac{\lambda\,v^\star}{2}\right|}{\sqrt{1 - \lambda + \frac{\lambda^2}{2}}} \ge \Phi^{-1}\left(1-\frac{\alpha}{2}\right) \right\} \\
& = 2\left[1 - \left(1-\frac{\alpha}{2}\right)\right] \\
& = \alpha.
\end{align*}
\end{proof}

It is important to note that although Proposition \ref{prop:TwoPointPC} establishes the validity and efficiency of using the mixture distribution as an approach to perform partial conditioning, by no means does Proposition \ref{prop:TwoPointPC} imply that using mixture distributions is the \emph{unique} approach to achieve partial conditioning.
}

\color{black}

\ifthenelse{1=1}{}{
\section*{Reviewer Comments}

Please motivate better the Hölder space-constrained estimation problem, e.g. in what modeling/practical situations does the constraint arise, and how does one interpret the parameters M and $\gamma$ so that the readers can see why and how this is a challenge.
-Second, the partial regression construction of the predictive IM in Section 3.2 reminds the reader of shrinkage estimation, an important connection that the authors should explore.
-Is the reader correct that as the paper stands, the proposed partial conditioning approach is restricted to handling the many means problem under constraints? If so, it’d be appropriate to reflect this in the title, since a general theory about partial conditioning for IM is not discussed.

In addition, the authors could improve their presentation of Section 4, either by streamlining some of the results or by adding additional comments to accompany them so the readers can grasp their significance. \textcolor{gray}{The following paragraph is furthermore duplicated and requires editing: ``Solving optimization problem (xxx) analytically requires finding the roots of a set of quadratic equations, a process which can become quite burdensome. Instead, we adopt numerical optimization methods such as BFGS to solve problem (xxx) in order to obtain the optimal mixing proportions which minimize the length of the plausibility interval.’'}

\subsection*{Reviewer #1}

Major comments
1. The terms "validity" and "efficiency" are emphasized in the abstract, Introduction, and Discussion without any indication that the two terms have different meanings than their standard definitions in the statistics literature.
2. The definition of "validity" is not connected with the standard meaning of that term in the statistics literature. Does it generalize the standard concept? Why satisfying that condition is important in practical data analysis is not explained.
3. Worse, "efficiency" is never formally defined.
4. IJAR readers would appreciate an explanation of the relation of the inferential models framework to the Dempster-Shaffer and imprecise probability models they are already familiar with.

Minor comments
1. In the figures, the axis labels and numbers on the ticks are too small.
2. The Introduction lacks a review of the strengths and limitations of the main previous approaches to the problem.
3. The Discussion lacks a clear presentation of the advantages and disadvantages of the proposed method compared to those previous approaches.

\subsection*{Reviewer #2}

IMs are not familiar to the general public, so, overall, Section 2 needs major improvements.
Section 2.1: Some further explanation on how to choose the predictive random sets in the P-step can be helpful, as it seems to be a key component of the IM construction.
An explanation of what (2.6) means, in words, can help the reader understand why the validity criteria is important. The presentation of the equivalent validity condition in terms of belief functions can also be beneficial.
Further explanation on why the consideration of lower and upper probabilities can be a good thing in statistical inference may be worth adding.
Section 2.2 on conditional IMs can be very hard to digest. Maybe further background on the dimension reduction advocated by Martin and Liu (2015a) would be beneficial. Moreover, further details and visual illustrations on the example (B=0) presented would benefit the reader.
I understand that the meaning and usefulness of marginal IMs can be read between the lines in the paper. But I believe some more direct explanation would benefit the reader. Perhaps a subsection?
Why do you say the chosen random set in the P-step on page 5 is "default"? This term was not used before…. If you are trying to mimic (2.2), shouldn't you also divide by ½ to have a standard normal? I understand this makes no difference in the plausibility/belief calculations below, but this can create some confusion…. As I said above, further explanation on how to choose a random set may be good.
\textcolor{gray}{The abbreviation "CIM" on page 5 can be deleted, as it was never used again.}
Conditionally admissible predictive random sets are introduced on page 5 with no explanation
Section 3.1: The discussion on why the fiducial solution may not work is interesting…. Can the authors add some explanation to the second paragraph?
Section 3.2: All of a sudden the authors use the term "predictive distribution". This should be defined to avoid confusion.
Section 3.2: A definition of a nested predictive random set would benefit the reader.

\textcolor{gray}{Not everybody will know what BFGS stands for on pages 11 and 13.
A better description of the simulation studies is needed, as well as the legends in some Figures. For example, what does Cond'(1pt) mean in Figure 3?
A better description of Figures 4, 5 and 6 is necessary. What are the "x" symbols in the graphs? What are the dotted lines?}
Section 5 deserves further explanations, such as the limitations of the proposed method, the main difficulties in implementing it, future directions, what a generalization of the proposed method can lead to and etc.

One more thing: In my opinion, the main advantage of the IM framework is its ability to provide valid probabilistic uncertainty quantification to any hypothesis about the unknown parameters, not just confidence intervals. Therefore, an interesting addition to the paper would be the comparison of IM's lower probabilities on some assertion of interest about the means with an alternative solution to the problem, e.g fiducial or Bayesian. How do they compare?
}

\section*{Appendix: Proofs}

\subsection*{Proof of Proposition \ref{prop:one-point-example}}

\begin{proof}
For the predictive random set $\mathcal{S} = [-|U_1|,|U_1|],\ U_1\sim\mathcal{N}(0,1)$, we have
$$ \mathsf{P}_\mathcal{S} \{\mathcal{S}\not\owns u^\star\} = \mathsf{P} \{|U_1|<|u_1^\star|\} = 2\Phi(|u_1^\star|)-1, $$
where $\Phi(\cdot)$ denotes the standard Normal CDF.
It follows that for each $\alpha\in(0,1)$,
$$ \mathsf{P}_{U_1} \{\mathsf{P}_\mathcal{S}\{\mathcal{S}\not\owns U_1\}\ge1-\alpha\} = \mathsf{P}_{U_1} \{2\Phi(|U_1|)-1\ge 1-\alpha\} = \mathsf{P}_{U_1} \{\Phi(|U_1|)\ge 1-\alpha/2\} = \alpha, $$
which verifies that \eqref{eq:PRSvalidity} holds. Thus, we have shown that the predictive random set $\mathcal{S}$ is valid, and Theorem~\ref{thm:validity} implies that the IM previously defined is valid as well.
\end{proof}

\subsection*{Proof of Proposition \ref{prop:3PointCond1}}

\begin{proof}
The conditional distributions for $U_1$ can be derived as follows:
\begin{align*}
\begin{dcases}
(U_1,V_{21})\sim\mathcal{N}_2 \left( \left[\begin{array}{c} 0 \\ 0 \end{array}\right], \left[\begin{array}{cc} 1 & -1 \\ -1 & 2 \end{array}\right] \right) & \quad \Rightarrow \quad U_1|V_{21}=v_{21} \sim \mathcal{N}\left(-\frac{v_{21}}{2}, \frac{1}{2}\right); \\
(U_1,V_{31})\sim\mathcal{N}_2 \left( \left[\begin{array}{c} 0 \\ 0 \end{array}\right], \left[\begin{array}{cc} 1 & -1 \\ -1 & 2 \end{array}\right] \right) & \quad \Rightarrow \quad U_1|V_{31}=v_{31} \sim \mathcal{N}\left(-\frac{v_{31}}{2}, \frac{1}{2}\right); \\
(U_1,V_{32})\sim\mathcal{N}_2 \left( \left[\begin{array}{c} 0 \\ 0 \end{array}\right], \left[\begin{array}{cc} 1 & 0 \\ 0 & 2 \end{array}\right] \right) & \quad \Rightarrow \quad U_1|V_{32}=v_{32} \sim \mathcal{N}(0,1).
\end{dcases}
\end{align*}
The conditional distributions for $U_2$ can be derived as follows:
\begin{align*}
\begin{dcases}
(U_2,V_{21})\sim\mathcal{N}_2 \left( \left[\begin{array}{c} 0 \\ 0 \end{array}\right], \left[\begin{array}{cc} 1 & 1 \\ 1 & 2 \end{array}\right] \right) & \quad \Rightarrow \quad U_2|V_{21}=v_{21} \sim \mathcal{N}\left(\frac{v_{21}}{2}, \frac{1}{2}\right); \\
(U_2,V_{31})\sim\mathcal{N}_2 \left( \left[\begin{array}{c} 0 \\ 0 \end{array}\right], \left[\begin{array}{cc} 1 & 0 \\ 0 & 2 \end{array}\right] \right) & \quad \Rightarrow \quad U_2|V_{31}=v_{31} \sim \mathcal{N}(0,1); \\
(U_2,V_{32})\sim\mathcal{N}_2 \left( \left[\begin{array}{c} 0 \\ 0 \end{array}\right], \left[\begin{array}{cc} 1 & -1 \\ -1 & 2 \end{array}\right] \right) & \quad \Rightarrow \quad U_2|V_{32}=v_{32} \sim \mathcal{N}\left(-\frac{v_{32}}{2}, \frac{1}{2}\right).
\end{dcases}
\end{align*}
The conditional distributions for $U_3$ can be derived as follows:
\begin{align*}
\begin{dcases}
(U_3,V_{21})\sim\mathcal{N}_2 \left( \left[\begin{array}{c} 0 \\ 0 \end{array}\right], \left[\begin{array}{cc} 1 & 0 \\ 0 & 2 \end{array}\right] \right) & \quad \Rightarrow \quad U_3|V_{21}=v_{21} \sim \mathcal{N}(0,1); \\
(U_3,V_{31})\sim\mathcal{N}_2 \left( \left[\begin{array}{c} 0 \\ 0 \end{array}\right], \left[\begin{array}{cc} 1 & 1 \\ 1 & 2 \end{array}\right] \right) & \quad \Rightarrow \quad U_3|V_{31}=v_{31} \sim \mathcal{N}\left(\frac{v_{31}}{2}, \frac{1}{2}\right); \\
(U_3,V_{32})\sim\mathcal{N}_2 \left( \left[\begin{array}{c} 0 \\ 0 \end{array}\right], \left[\begin{array}{cc} 1 & 1 \\ 1 & 2 \end{array}\right] \right) & \quad \Rightarrow \quad U_3|V_{32}=v_{32} \sim \mathcal{N}\left(\frac{v_{32}}{2}, \frac{1}{2}\right).
\end{dcases}
\end{align*}
The proof is complete by recalling that $v_{ij} = u_i-u_j$ for $1\le j<i\le 3$.
\end{proof}

\subsection*{Proof of Corollary \ref{cor:3PointCond1}}

\begin{proof}
We notice that the results in \eqref{eq:3PointBounds} imply the following:
\begin{align*}
u_1-u_2 = -v_{21} & \in [Y_1-Y_2-B_{21}, Y_1-Y_2+B_{21}], \\
u_1-u_3 = -v_{31} & \in [Y_1-Y_3-B_{31}, Y_1-Y_3+B_{31}]; \\
u_2-u_1 = v_{21} & \in [Y_2-Y_1-B_{21}, Y_2-Y_1+B_{21}], \\
u_2-u_3 = -v_{32} & \in [Y_2-Y_3-B_{32}, Y_2-Y_3+B_{32}]; \\
u_3-u_1 = v_{31} & \in [Y_3-Y_1-B_{31}, Y_3-Y_1+B_{31}], \\
u_3-u_2 = v_{32} & \in [Y_3-Y_2-B_{32}, Y_3-Y_2+B_{32}].
\end{align*}
Combining these results with Proposition \ref{prop:3PointCond1} completes the proof of the corollary.
\end{proof}

\subsection*{Proof of Proposition \ref{prop:3PointCond2} }

\begin{proof}
The conditional distributions for $U_1$ can be derived as follows:
\begin{align*}
\begin{dcases}
(U_1,V_{21},V_{31}) \sim \mathcal{N}_3 \left( \left[\begin{array}{c} 0 \\ 0 \\ 0 \end{array}\right], \left[\begin{array}{ccc} 1 & -1 & -1 \\ -1 & 2 & 1 \\ -1 & 1 & 2 \end{array}\right] \right) & \quad \Rightarrow \quad U_1|v_{21},v_{31} \sim \mathcal{N}\left(-\dfrac{v_{21}+v_{31}}{3},\dfrac{1}{3}\right); \\
(U_1,V_{21},V_{32}) \sim \mathcal{N}_3 \left( \left[\begin{array}{c} 0 \\ 0 \\ 0 \end{array}\right], \left[\begin{array}{ccc} 1 & -1 & 0 \\ -1 & 2 & -1 \\ 0 & -1 & 2 \end{array}\right] \right) & \quad \Rightarrow \quad U_1|v_{21},v_{32} \sim \mathcal{N}\left(-\dfrac{2v_{21}+v_{32}}{3},\dfrac{1}{3}\right); \\
(U_1,V_{31},V_{32}) \sim \mathcal{N}_3 \left( \left[\begin{array}{c} 0 \\ 0 \\ 0 \end{array}\right], \left[\begin{array}{ccc} 1 & -1 & 0 \\ -1 & 2 & 1 \\ 0 & 1 & 2 \end{array}\right] \right) & \quad \Rightarrow \quad U_1|v_{31},v_{32} \sim \mathcal{N}\left(\dfrac{v_{32}-2v_{31}}{3},\dfrac{1}{3}\right).
\end{dcases}
\end{align*}
The conditional distributions for $U_2$ can be derived as follows:
\begin{align*}
\begin{dcases}
(U_2,V_{21},V_{31}) \sim \mathcal{N}_3 \left( \left[\begin{array}{c} 0 \\ 0 \\ 0 \end{array}\right], \left[\begin{array}{ccc} 1 & 1 & 0 \\ 1 & 2 & 1 \\ 0 & 1 & 2 \end{array}\right] \right) & \quad \Rightarrow \quad U_2|v_{21},v_{31} \sim \mathcal{N}\left(\dfrac{2v_{21}-v_{31}}{3},\dfrac{1}{3}\right); \\
(U_2,V_{21},V_{32}) \sim \mathcal{N}_3 \left( \left[\begin{array}{c} 0 \\ 0 \\ 0 \end{array}\right], \left[\begin{array}{ccc} 1 & 1 & -1 \\ 1 & 2 & -1 \\ -1 & -1 & 2 \end{array}\right] \right) & \quad \Rightarrow \quad U_2|v_{21},v_{32} \sim \mathcal{N}\left(\dfrac{v_{21}-v_{32}}{3},\dfrac{1}{3}\right); \\
(U_2,V_{31},V_{32}) \sim \mathcal{N}_3 \left( \left[\begin{array}{c} 0 \\ 0 \\ 0 \end{array}\right], \left[\begin{array}{ccc} 1 & 0 & -1 \\ 0 & 2 & 1 \\ -1 & 1 & 2 \end{array}\right] \right) & \quad \Rightarrow \quad U_2|v_{31},v_{32} \sim \mathcal{N}\left(\dfrac{v_{31}-2v_{32}}{3},\dfrac{1}{3}\right).
\end{dcases}
\end{align*}
The conditional distributions for $U_3$ can be derived as follows:
\begin{align*}
\begin{dcases}
(U_3,V_{21},V_{31}) \sim \mathcal{N}_3 \left( \left[\begin{array}{c} 0 \\ 0 \\ 0 \end{array}\right], \left[\begin{array}{ccc} 1 & 0 & 1 \\ 0 & 2 & 1 \\ 1 & 1 & 2 \end{array}\right] \right) & \quad \Rightarrow \quad U_3|v_{21},v_{31} \sim \mathcal{N}\left(\dfrac{2v_{31}-v_{21}}{3},\dfrac{1}{3}\right); \\
(U_3,V_{21},V_{32}) \sim \mathcal{N}_3 \left( \left[\begin{array}{c} 0 \\ 0 \\ 0 \end{array}\right], \left[\begin{array}{ccc} 1 & 0 & 1 \\ 0 & 2 & -1 \\ 1 & -1 & 2 \end{array}\right] \right) & \quad \Rightarrow \quad U_3|v_{21},v_{32} \sim \mathcal{N}\left(\dfrac{v_{21}+2v_{32}}{3},\dfrac{1}{3}\right); \\
(U_3,V_{31},V_{32}) \sim \mathcal{N}_3 \left( \left[\begin{array}{c} 0 \\ 0 \\ 0 \end{array}\right], \left[\begin{array}{ccc} 1 & 1 & 1 \\ 1 & 2 & 1 \\ 1 & 1 & 2 \end{array}\right] \right) & \quad \Rightarrow \quad U_3|v_{31},v_{32} \sim \mathcal{N}\left(\dfrac{v_{31}+v_{32}}{3},\dfrac{1}{3}\right).
\end{dcases}
\end{align*}
It can be easily verified that the three conditional distributions in each group are actually equivalent, and that the results \eqref{eq:3PointCond1} -- \eqref{eq:3PointCond3} hold.
\end{proof}

\subsection*{Proof of Corollary \ref{cor:3PointCond2} }

\begin{proof}
We notice that the results in \eqref{eq:3PointBounds} implies the following:
\begin{align*}
2u_1-u_2-u_3 = -(v_{21}+v_{31}) & \in [2Y_1-Y_2-Y_3-(B_{21}+B_{31}), 2Y_1-Y_2-Y_3+(B_{21}+B_{31})]; \\
2u_2-u_1-u_3 = v_{21}-v_{32} & \in [2Y_2-Y_1-Y_3-(B_{21}+B_{32}), 2Y_2-Y_1-Y_3+(B_{21}+B_{32})]; \\
2u_3-u_1-u_2 = v_{31}+v_{32} & \in [2Y_3-Y_1-Y_2-(B_{31}+B_{32}), 2Y_1-Y_2-Y_3+(B_{31}+B_{32})].
\end{align*}
Combining these results with Proposition \ref{prop:3PointCond2} completes the proof of the corollary.
\end{proof}

\ifthenelse{1=0}{}{
  \bibliographystyle{plain}
  \bibliography{nonpar}
}

\end{document}